\documentclass[a4paper, 11pt]{amsart}
\usepackage{amssymb,amsthm,graphicx,accents,xcolor,hyperref}
\usepackage{tikz}
\usepackage[utf8]{inputenc}
\usepackage[T1]{fontenc}

\theoremstyle{plain}
\newtheorem{theorem}{Theorem}[section]
\newtheorem{lemma}[theorem]{Lemma}
\newtheorem{proposition}[theorem]{Proposition}
\newtheorem{corollary}[theorem]{Corollary}
\newtheorem{question}[theorem]{Question}

\theoremstyle{definition}
\newtheorem{definition}[theorem]{Definition}
\newtheorem{example}[theorem]{Example}

\theoremstyle{remark}
\newtheorem{remark}[theorem]{Remark}

\newtheorem*{proof*}{Proof}

\author{Bruno Teheux}
\address{Department of Mathematics, University of Luxembourg, Maison du
Nombre, 6 avenue de la Fonte, L-4364 Esch-sur-Alzette, Luxembourg}
\email{bruno.teheux@uni.lu}
\title[Dualities and completions of median algebras]{Dualities and
completions of median algebras and semilattices: a universal-algebraic
approach}
\date{9 November 2026}

\newcommand{\Z}{\mathbb{Z}}
\newcommand{\N}{\mathbb{N}}
\newcommand{\R}{\mathbb{R}}
\newcommand{\utilde}[1]{\smash[b]{\underaccent{\tilde}{#1}}}
\renewcommand{\emptyset}{\varnothing}
\newcommand{\co}{{\bullet}}
\newcommand{\alg}[1]{\mathbf{#1}}
\newcommand{\cate}[1]{\mathsf{#1}}
\newcommand{\Poc}{\cate{Poc}}
\newcommand{\M}{\cate{M}}
\newcommand{\classop}[1]{\mathbb{#1}}

\newcommand{\ST}{\mathrm{ST}}
\newcommand{\Trans}{\mathrm{Trans}}
\newcommand{\Mod}{\mathrm{Mod}}
\newcommand{\Fr}{\mathrm{Fr}}
\newcommand{\clK}{\mathcal{K}}
\DeclareMathOperator{\rank}{rank}

\let\svhat\hat
\renewcommand{\hat}[1]{\smash{\svhat{#1}}}
\let\svoverline\overline
\renewcommand{\overline}[1]{\smash{\svoverline{#1}}}

\begin{document}
\begin{abstract}
We show that universal algebra, and the theory of natural dualities in
particular, provides a uniform framework for the dualities between median
algebras and bounded poc sets and their Boolean topological variants. We
identify the Roller compactification of a $\mathrm{CAT}(0)$ cube complex with
the profinite completion and the natural extension of the median algebra of
its vertices. We then develop the order-theoretic perspective, which leads to
a characterization of the profinite completion by the property that its
intervals are the canonical extensions of the original ones. Since the
basepoint orders of a median algebra have no top in general, we introduce
conditional canonical extensions, show that they are unique, and realize
them inside the profinite completion. Finally, a transfer theorem carries
universal sentences from the finite quotients of a median algebra to its
profinite completion, and we close with open questions on a correspondence
theory between median algebras and their bounded poc sets.
\end{abstract}
\maketitle
\tableofcontents
\section{Introduction}

Median algebras and their order-theoretic counterparts, the median
semilattices, were introduced in lattice theory in the 1940s and 1950s. Their
duals, the bounded poc sets, appeared around 1980 together with the first
dualities. All three resurfaced in the 1990s in geometric group theory, which
found them in trees and in $\mathrm{CAT}(0)$ cube complexes.
Section~\ref{subsec:history2} recalls this history.

The traditions have developed, each in its own vocabulary, the same two kinds of
tools. The first one is duality between median algebras and bounded poc sets.
The dualities of Isbell \cite{Isbell1980} and Werner \cite{Werner1981} reappear
in the work of Roller \cite{Roller2016} and Bowditch \cite{Bowditch2022}.  The
second one is the extension of a median algebra by the points that it misses.
The Roller compactification of a $\mathrm{CAT}(0)$ cube complex is, on the
algebraic side, the natural extension and the profinite completion of the median
algebra of its vertices, and the hyperplane-collapses of Genevois
\cite[\S2.3]{Genevois2024}, or restriction quotients of Caprace and Sageev
\cite[\S2.3]{CapraceSageev2011}, are the finite quotients from which the latter
is built. The
zero-completion of Bandelt and Meletiou \cite{BandeltMeletiou1993}, in the
interval-based form given to it by Fioravanti \cite{Fioravanti2020} and Bowditch
\cite{Bowditch2022}, sits next to the canonical extensions of lattice theory.
These constructions are usually studied separately, and the bridges between them
are the subject of this paper.

This paper has two aims. The first one is to show that universal algebra
provides a unified setting for constructions that the geometric literature
meets one at a time. A few general theorems, applied to the two element
median algebra $\alg{2}$, deliver all of them at once, and it is the theory
of natural dualities \cite{ClarkDavey} that supplies most of these theorems.
The second aim is to characterize these constructions. We follow the three
points of view on order
completions described by Fussner and Gehrke \cite{FussnerGehrke2025}, namely
a concrete construction, an internal characterization saying how the original
structure sits in the completion, and a universal property presenting the
completion as a left adjoint. Throughout, median semilattices are treated on
a par with median algebras, and each completion is also studied as an order
completion of the semilattices $(\alg{A}, \leq_a)$, a point of view that, as
far as we know, the geometric literature does not develop.

Section~\ref{sec:dualities} obtains the two dualities by turning the crank of
natural duality theory. Since the median operation is a majority term, a
general theorem of the theory yields a strong duality between median algebras
and Boolean topological bounded poc sets, and the TopSwap
Theorem \cite{DaveyHaviarPriestley2012} turns it into a dual
equivalence between Boolean topological median algebras and bounded poc sets
(Theorem~\ref{thm:dualities}). Neither duality is new. Here both come out of
the general theory, the second one as a formal consequence of the first.
The Boolean topological median algebras are the profinite ones
(Proposition~\ref{prop:profinite-char}), by the theory of standard
topological algebras \cite{ClarkDaveyFreeseJackson2004}.\footnote{The class
$\classop{IS}_c\classop{P}^+(\alg{M}_\tau)$ of Boolean topological algebras
generated by a finite algebra $\alg{M}$ is \emph{standard} if it consists of
all the Boolean topological algebras that satisfy the universal Horn theory of
$\alg{M}$ \cite{ClarkDaveyFreeseJackson2004}.}

Section~\ref{sec:profinite} studies the profinite completion $\hat{\alg{A}}$
of a median algebra $\alg{A}$. Two further general theorems identify it with
the natural extension of $\alg{A}$ and with the algebra $(\alg{A}^+)_+$ of
morphisms from the bounded poc set $\alg{A}^+$ of prime convex subsets of
$\alg{A}$ into the two element bounded poc set
(Proposition~\ref{prop:natext}, Theorem~\ref{thm:DGHP} and
Corollary~\ref{prop:profinite}). Read through the
vocabulary of geometric group theory, this says that the Roller
compactification of a $\mathrm{CAT}(0)$ cube complex is the profinite
completion of the median algebra of its vertices, that is, the limit of its
finite hyperplane-collapses (Subsection~\ref{subsec:roller}). We have not
found this identification in either literature. It comes with
functoriality, the profinite completion being left adjoint to the functor
forgetting the topology (Proposition~\ref{prop:ump-reflection}). The duality
also describes all the Boolean topological median algebras in which $\alg{A}$
has a dense image,
through the bounded sub-poc sets of $\alg{A}^+$, and shows that the Roller
compactification is the only one in which $\alg{A}$ embeds densely when the
intervals of $\alg{A}$ are finite (Proposition~\ref{prop:compactifications}).

Section~\ref{sec:can-ext} gives the internal characterization. Each interval
of $\hat{\alg{A}}$ is the canonical extension \cite{Dunn2005} of the
corresponding interval of $\alg{A}$
(Theorem~\ref{thm:can-ext}), and this property singles out the profinite
completion among the Boolean topological median algebras in which $\alg{A}$
is dense (Theorem~\ref{thm:char-profinite}). The zero-completion embeds into
$\hat{\alg{A}}$ and coincides with it when the intervals of $\alg{A}$ are
finite. In general it is smaller, since it leaves the intervals of $\alg{A}$
untouched whereas $\hat{\alg{A}}$ completes each of them
(Remark~\ref{rem:ccan-vs-zero}).

Section~\ref{sec:ccan} turns to the whole median semilattice $(\alg{A},
\leq_a)$, which has in general no top. The order $(\hat{\alg{A}}, \leq_a)$ is
then in general not a canonical extension of it, and it may even fail to be
dense (Example~\ref{ex:escalator}), as can the profinite completions of
semilattices \cite{Gouveia2013}. We introduce conditional canonical
extensions, the dense and compact conditional completions in which the
original poset is cofinal, prove that they are unique
(Theorem~\ref{thm:unique-ccan}), and find the one of $(\alg{A}, \leq_a)$
inside $\hat{\alg{A}}$, as the union of the intervals $[a;b]$ with $b \in A$
(Theorem~\ref{thm:realization}). This union is the whole of $\hat{\alg{A}}$
exactly when $\alg{A}$ is cofinal in $(\hat{\alg{A}}, \leq_a)$. The
conditional canonical extension is left adjoint to a forgetful functor
(Corollary~\ref{cor:ccan-adjunction}), which extends the description of the
canonical extension of a bounded distributive lattice as the free completely
distributive lattice over it.

Section~\ref{sec:transfer} is a first application. Every universal sentence
that holds in the finite quotients of $\alg{A}$ holds in $\hat{\alg{A}}$ and
in all its subalgebras (Theorem~\ref{thm:transfer}), and the non-strict
rank bounds of the geometric literature on the Roller compactification and
on the components of its boundary become instances of this transfer
(Corollary~\ref{cor:rank-transfer}).
Section~\ref{sec:conclusion} closes with research perspectives, among which a
correspondence theory relating sentences about median algebras to
properties of their bounded poc sets.

\section{Median algebras, poc sets and median semilattices}\label{subsec:history2}

The three structures that we investigate in this paper are median algebras,
median semilattices, which are their basepointed form, and bounded poc sets.
All three are classical. Median algebras and bounded poc sets descend from two
independent traditions, whereas median semilattices belong to the first one.

The first line of descent runs through lattice theory. Birkhoff and Kiss
\cite{BirkhoffKiss1947} started from the ternary term $(a\wedge b) \vee (a
\wedge c) \vee (b \wedge c)$ of a lattice and observed that it obeys the median
laws exactly when the lattice is distributive, so that the variety of median
algebras and that of distributive lattices share the same median equational
theory. Independently, and from betweenness rather than from lattices,
Sholander axiomatized the same class \cite[\S8]{Sholander1952} and went on to
study it in \cite{Sholander1954, Sholander1954b}. Grau \cite{Grau1947}
axiomatized the Boolean case, his \emph{ternary Boolean algebras} carrying,
besides the median, a unary complement subject to $m(a,b,b^\co) = a$; he showed
them to correspond exactly to Boolean algebras. The representation of median
algebras by systems of prime convex subsets was then developed by Isbell
\cite{Isbell1980}, Bandelt and Hedl\'{\i}kov\'a \cite{BandeltHedlikova1983} and
van de Vel \cite{vandeVel1993}, well before the notion was taken up in geometric
group theory.

Median algebras and poc sets resurfaced in geometric group theory, which
brought to them new examples, a new motivation and a distinct vocabulary
\cite{Roller2016, Bowditch2022}. Remark~\ref{rem:terminology} below matches the
two vocabularies.

The present section fixes the three notions and the constructions relating
them.

\subsection{Median algebras}\label{subsec:medpoc}

\begin{definition}\label{def:median}
A \emph{median algebra} is a ternary algebra $\alg{A} = (A,m)$ in which the
\emph{symmetry}, \emph{majority} and \emph{distributivity} axioms
\begin{gather}
m(a,b,c) = m(b,a,c) = m(b,c,a)\label{eq:symmetry}\\
m(a,a,b) = a\label{eq:majority}\\
m(m(a,b,c), d, e) = m(a, m(b, d, e), m(c, d, e))\label{eq:distributivity}
\end{gather}
hold. We denote by $\M$ the variety of median algebras defined by these
equations.
\end{definition}

The three axioms are, verbatim, the postulates $(3.3)$, $(3.2)$ and $(3.4)$ of
Birkhoff and Kiss \cite[p.~751]{BirkhoffKiss1947}, who add one further
postulate, $(0,a,I) = a$ for a fixed pair of constants, in order to pin down
bounded distributive lattices rather than median algebras. Sholander
axiomatized the same class by the two identities $m(a,a,b) = a$ and
$m(m(a,b,c), m(a,b,d), e) = m(m(c,d,e),a,b)$ \cite[\S8]{Sholander1952}. The
equational base above is the one adopted in the surveys of Isbell
\cite{Isbell1980} and of Bandelt and Hedl\'{\i}kov\'a
\cite{BandeltHedlikova1983}.

We use two basic examples. The first one is the algebra $\alg{2}$
obtained from $\{0,1\}$ by adding the majority operation, that is, the
operation $m$ determined by $m(x,x,y) = x$. The second one is a tree $T$
with vertex set $V$. Any two vertices $a$
and $b$ of $T$ are joined by a unique geodesic $[a,b]$, and for any three
vertices $a$, $b$ and $c$ the intersection $[a,b] \cap [a,c] \cap [b,c]$
consists of a single vertex, denoted by $m(a,b,c)$. The algebra $\alg{T} =
(V,m)$ obtained in this way is a median algebra.
Sholander axiomatized trees in this language, and characterized them among
median algebras by one further postulate \cite[\S2 and \S9]{Sholander1952}.

\begin{definition}\label{def:convex}
A subset $x$ of a median algebra $\alg{A}$ is called \emph{convex} if
\begin{equation}\label{eq:conv}
    \text{at least two elements among } a,b,c \text{ are in } x \quad \implies
\quad m(a,b,c) \in x,
\end{equation}
and \emph{prime convex} if it is convex and satisfies the converse of
\eqref{eq:conv}.
\end{definition}

Both $\emptyset$ and $A$ are prime convex subsets of any median algebra
$\alg{A}$, and a convex subset $x$ of $\alg{A}$ is prime if and only if $A
\setminus x$ is convex. For $a, b \in A$, the \emph{interval} $[a;b]$ is the
set $\{x \in A \mid m(a,x,b) = x\}$. A median algebra is \emph{discrete} if
all its intervals are finite \cite[\S11.1]{Bowditch2022},
\cite[pp.~23--24]{Roller2016}.

We often identify without further notice a prime convex subset $x$ of a median
algebra $\alg{A}$ with its characteristic function, which is a homomorphism $x
\colon \alg{A} \to \alg{2}$. This identification is a bijection, whose inverse
associates the prime convex subset $x^{-1}(1)$ to any $x \in \M(\alg{A},
\alg{2})$. The following lemma is the median counterpart of the prime ideal
theorem.

\begin{lemma}[Extension Lemma (AC)]\label{lem:sep}
Let $\alg{A}$ be a median algebra, let $x$ be a convex subset of $\alg{A}$ and
let $b$ be an element of $A \setminus x$. There is a prime convex subset of
$\alg{A}$ that contains $x$ and avoids $b$.
\end{lemma}
\begin{proof}
Zorn's Lemma provides a convex subset of $\alg{A}$ that contains $x$, avoids
$b$ and is maximal with respect to these two properties. Maximality forces that
subset to be prime.
\end{proof}

The lemma goes back to Nieminen \cite[Theorem~1]{Nieminen1978}, who derived it
from a separation theorem for semilattices of Balbes
\cite[Theorem~2.2]{Balbes1969}; see also
\cite[Lemma~1.4]{BandeltHedlikova1983} and \cite[Theorem~2.7]{Roller2016}.
Throughout, we use the class operators of universal algebra, writing
$\classop{I}$, $\classop{H}$, $\classop{S}$ and $\classop{P}$ for the operators
forming, from a class, its isomorphic copies, its homomorphic images, its
substructures and its products, and adding a superscript $+$ to $\classop{P}$
to restrict to products over nonempty index sets.

Applied to the convex subset $\{a\}$ and to an element $b \neq a$,
Lemma~\ref{lem:sep} shows
that the homomorphisms from a median algebra into $\alg{2}$ separate its
points, whence the following representation.

\begin{theorem}[AC]\label{prop:med-axioms}
The variety of median algebras is the class $\classop{ISP}(\alg{2})$. More
precisely, every median algebra $\alg{A}$ embeds into the power
$\alg{2}^{\M(\alg{A}, \alg{2})}$ under the evaluation map sending $a \in A$ to
the family $(h(a))_{h \in \M(\alg{A}, \alg{2})}$.
\end{theorem}

The form $\M = \classop{ISP}(\alg{2})$ of this classical representation
\cite[Theorem~1.5]{BandeltHedlikova1983} is the starting point of natural
duality theory (Section~\ref{sec:dualities}). In particular, the variety
$\M$ is \emph{locally finite}, that is, its finitely generated members are
finite. Indeed, an element of the subalgebra of a power of $\alg{2}$ generated
by $n$ elements is the value of a term at the generators, computed
coordinatewise, so it depends only on the function from $\{0,1\}^n$ to
$\{0,1\}$ that this term induces on $\alg{2}$. Hence a median algebra
generated by $n$ elements has at most $2^{2^n}$ elements
\cite[Proposition~3.3.3]{Bowditch2022}.

\subsection{Bounded poc sets}\label{subsec:poc}

Deleting an edge of a tree $T$ splits its vertex set into two prime convex
subsets of $\alg{T}$, called \emph{half-trees}, and we write
$\mathcal{H}(T)$ for the set of half-trees of $T$ together with $\emptyset$ and
$V$, ordered by inclusion. The
complement in $V$ of a half-tree is again a half-tree, so that complementation
is an involution of $\mathcal{H}(T)$ that reverses inclusion and exchanges
$\emptyset$ and $V$. Moreover, no member of $\mathcal{H}(T)$ other than
$\emptyset$ is contained in its own complement. The next definition turns these
properties into axioms.

\begin{definition}\label{def:poc}
A \emph{poc set} $(X,\leq, {}^\co)$ is given by a partially ordered set $(X, \leq)$
equipped with a unary involutive order reversing operation ${}^\co$ that
satisfies $x \not\leq x^\co$ for every $x\in X$.

A \emph{bounded poc set} $(X, \leq, {}^\co, 0, 1)$ is given by a bounded poset
$(X,\leq, 0,1)$ with $0 \neq 1$, equipped with a unary involutive order
reversing operation ${}^\co$ that satisfies $0^\co =1$, and $x \leq x^\co$
if and only if $x=0$. We denote
by $\Poc$ the category of bounded poc sets with structure preserving maps as
arrows.
\end{definition}

Bounded poc sets are Isbell's \emph{binary messages} \cite[\S6]{Isbell1980},
whose definition allows in addition the one element structure, and they also
appear in Werner \cite{Werner1981}. The name \emph{poc set} and the axioms in
the form used here are Roller's \cite[\S1.3]{Roller2016}; his poc sets carry a
smallest element, and are therefore the bounded ones.

Deleting the bounds of a bounded poc set gives a poc set, and every poc set
arises in this way from a bounded poc set, unique up to isomorphism. This
correspondence does not extend to morphisms, since a morphism of bounded poc
sets may send elements other than the bounds to bounds, as the morphisms into
the two element bounded poc set do.

We denote by $\utilde{\alg{2}}$ the two element bounded poc set. The category
$\Poc$ is generated by $\utilde{\alg{2}}$ in the following sense.

\begin{proposition}[AC]\label{prop:char}
Every bounded poc set is isomorphic to a substructure of a power of
$\utilde{\alg{2}}$ over a nonempty index set, that is, $\Poc =
\classop{ISP}^+(\utilde{\alg{2}})$.
\end{proposition}

The proposition\footnote{The empty index set has to be excluded, since the
empty product is the one element structure, in which $0$ and $1$ coincide.}
is Isbell's Theorem~6.10 \cite{Isbell1980}, and its proof is
yet another application of Zorn's Lemma.

We now pass from bounded poc sets to median algebras. The majority operation
of $\{0,1\}$ is idempotent, monotone in each argument and self-dual, in the
sense that $m(1-p,1-q,1-r) = 1-m(p,q,r)$, so the pointwise median of three
morphisms $\alg{X} \to \utilde{\alg{2}}$ is again
a morphism \cite[proof of Theorem~6.12]{Isbell1980}.

\begin{definition}\label{def:complex}
The \emph{complex algebra} of a bounded poc set $\alg{X}$ is the subalgebra
$\alg{X}_+ = (\Poc(\alg{X}, \utilde{\alg{2}}), m)$ of $\alg{2}^X$, whose median
is computed pointwise.
\end{definition}

A morphism $a \colon \alg{X} \to \utilde{\alg{2}}$ is the characteristic
function of $a^{-1}(1)$, which is an \emph{ultrafilter} of $\alg{X}$, that is,
a nonempty proper upset of $(X, \leq)$ that contains exactly one of $x$ and
$x^\co$ for every $x \in X$. Every ultrafilter arises in this way, and we
identify without further notice each morphism with its ultrafilter. The median
of $\alg{X}_+$ then reads
\[
  m(a,b,c) = (a\cap b) \cup (a \cap c) \cup (b\cap c).
\]
In geometric group theory an ultrafilter is called a \emph{flow}
\cite[\S9]{Bowditch2022}. Roller keeps the word ultrafilter for the notion
above and reserves \emph{orientation} for the weaker one obtained by dropping
the requirement that the set be an upset \cite[\S3]{Roller2016}.

For ultrafilters $a$, $b$ and $c$ of $\alg{X}$, the description of the median
just given turns the condition $c \in [a;b]$, computed in $\alg{X}_+$, into
\[
  a \cap b \;\subseteq\; c \;\subseteq\; a \cup b.
\]

In the other direction, the prime convex subsets of a median algebra form a
bounded poc set.

\begin{lemma}[{\cite[proof of Theorem~6.12]{Isbell1980}}]\label{lem:canonical-poc}
If $\alg{A}$ is a median algebra, then the set of prime convex subsets of
$\alg{A}$, ordered by inclusion, with the set complement in $A$ as involution
and with $\emptyset$ and $A$ as bounds, is a bounded poc set, called \emph{the
canonical bounded poc set of $\alg{A}$} and
denoted by
$\alg{A}^+$.
\end{lemma}
\begin{proof}
The verification is again routine; see also \cite[\S4]{Roller2016}.
\end{proof}

The canonical bounded poc set $\alg{A}^+$ is a substructure of
$\utilde{\alg{2}}^A$.

For a tree $T$, the canonical bounded poc set $\alg{T}^+$ is the half-tree
system $\mathcal{H}(T)$, and $(\alg{T}^+)_+$ consists of the evaluations at the
vertices of $T$ and of its ends \cite[Examples~4.3(vi)]{Roller2016}. For an
arbitrary median algebra $\alg{A}$, Corollary~\ref{prop:profinite} identifies
$(\alg{A}^+)_+$ with the profinite completion of $\alg{A}$.

\begin{remark}[Terminology in geometric group theory]\label{rem:terminology}
The terminology of the geometric tradition differs from ours. A convex subset
whose complement is again convex, a prime convex subset in our terminology, is
what the theory of median algebras and of abstract convexity calls a
\emph{halfspace}, and the canonical bounded poc set $\alg{A}^+$ of a median
algebra $\alg{A}$ is accordingly called its \emph{halfspace system}, whose
elements are the halfspaces of $\alg{A}$, ordered by inclusion, with set
complement as involution (see, \emph{e.g.}, \cite[\S2 and
Convention~4.2]{Roller2016}). Symmetrically, the complex algebra $\alg{X}_+$
of a bounded poc set $\alg{X}$ is the median algebra of its ultrafilters.
When $\alg{X}$ arises as the halfspace system of a \emph{space
with walls} in the sense of Haglund and Paulin \cite{HaglundPaulin1998} and
Nica \cite{Nica2004}, see also \cite[\S19]{Bowditch2022}, this construction
produces the vertices of the associated $\mathrm{CAT}(0)$ cube complex via
Sageev's construction \cite{Sageev1995}.
\end{remark}

\subsection{Median semilattices}

For any $\alg{A} \in \M$ and any $a \in A$, define the relation $\leq_a$ on $A$ by
\[
  b \leq_a c \quad\text{if and only if}\quad m(a,b,c) = b.
\]
This is a partial order on $A$ with bottom element $a$, called the \emph{base
point}. The binary operation
\[
  b \wedge_a c := m(a,b,c)
\]
makes $(\alg{A}, \leq_a)$ a meet-semilattice. We use the notations
$(\alg{A}, \leq_a)$ and $(\alg{A}, \wedge_a)$ interchangeably.
This meet-semilattice satisfies:
\begin{enumerate}
\item\label{item:mediansemi01} every principal ideal of $(\alg{A}, \leq_a)$
  is a distributive lattice;
\item\label{item:mediansemi02} any three elements of $A$ admit a join in
  $(\alg{A}, \leq_a)$ whenever each pair of them has an upper bound.
\end{enumerate}

\begin{lemma}\label{lem:interval-dual}
Let $\alg{A} \in \M$ and let $a, b, c, d \in A$ with $b, c \in [a;d]$.
Then $b \leq_a c$ if and only if $c \leq_d b$; that is, the posets
$([a;d], \leq_a)$ and $([a;d], \leq_d)$ are
order-duals of each other.
\end{lemma}
\begin{proof}
The distributivity axiom for $m$, together with the hypothesis $c \in [a;d]$,
gives
\[
  m\bigl(m(a,b,c),\, c,\, d\bigr)
  = m\bigl(c,\, m(a,c,d),\, m(b,c,d)\bigr)
  = m\bigl(c,\, c,\, m(b,c,d)\bigr)
  = c.
\]
If $b \leq_a c$, meaning $m(a,b,c) = b$, the left-hand side equals $m(b,c,d)$.
Hence $m(b,c,d) = c$, which means that $c \leq_d b$. The converse follows
by symmetry.
\end{proof}

\begin{definition}\label{def:median-semilattice}
A \emph{median semilattice} is a meet-semilattice with bottom element in which
every principal ideal is a distributive lattice and any three elements admit a
join whenever each pair of them has an upper bound.
\end{definition}

Median semilattices were introduced by Sholander
\cite{Sholander1954b} (see also \cite[\S5]{BandeltHedlikova1983}). Conversely,
any median semilattice
$(A, \wedge)$ yields a median algebra on $A$ by setting
\begin{equation}\label{eq:semi-to-med}
  m(a,b,c) := (a \wedge b) \vee (a \wedge c) \vee (b \wedge c),
\end{equation}
distributive lattices being the special case with a top element. The median
algebra recovered from $(\alg{A}, \wedge_a)$ via \eqref{eq:semi-to-med}
coincides with $\alg{A}$, for any
$a \in A$ \cite[Theorem~(3.1)]{Sholander1954b}, see also \cite{Evans1982}.

\begin{lemma}\label{lem:discrete-ideals}
Let $\alg{A} \in \M$ and $a \in A$. Then $\alg{A}$ is discrete if and only if
every principal ideal of $(\alg{A}, \leq_a)$ is finite.
\end{lemma}
In the language of walls this is the transitivity of almost equality
\cite[\S8]{Roller2016}, and we give the short order-theoretic proof.
\begin{proof}
The principal ideals of $(\alg{A}, \leq_a)$ are the intervals $[a;b]$,
which settles the forward direction. Conversely, let $b, c \in A$. By
\eqref{eq:semi-to-med}, every $x \in [b;c]$ satisfies $x = m(b,c,x) = (x
\wedge_a b) \vee (x \wedge_a c) \vee (b \wedge_a c)$, and the last term may be
dropped because $b \wedge_a c \leq_a x \wedge_a b$. Hence $x$ is determined by
the pair $(x \wedge_a b, x \wedge_a c)$ of ${\downarrow}b \times
{\downarrow}c$, and $[b;c]$ is finite as soon as ${\downarrow}b$ and
${\downarrow}c$ are.
\end{proof}

\section{Categorical dualities}\label{sec:dualities}

The constructions $-^+\colon \cate{M} \to \cate{Poc}$ and $-_+\colon \cate{Poc}
\to \cate{M}$ have so far been defined on objects only, and are therefore not
yet functorial. It is possible to lift them to functors in a way that leads to
categorical dualities, and the lifting is not unique. It gives rise to two
distinct categorical dualities, one between $\cate{M}$ and the category
$\Poc_\tau$ of Boolean topological bounded poc sets, and a partner duality
between the category $\cate{M}_\tau$ of Boolean topological median algebras and
$\cate{Poc}$. They are the two items of Theorem~\ref{thm:dualities}.

Neither duality is new. The first is due to Werner \cite{Werner1981} and,
independently, to Isbell \cite{Isbell1980}, in both cases as a special case of
a more general duality, and Clark and Davey state it as a strong duality
\cite[Theorem~4.3.4]{ClarkDavey}. What this section adds is that the second
one is a formal consequence of the first.

Observe that $\cate{M} = \classop{ISP}(\alg{2})$ and
$\cate{Poc}=\classop{ISP}^+(\utilde{\alg{2}})$ by
Theorem~\ref{prop:med-axioms} and Proposition~\ref{prop:char}, and that $-^+$ and $-_+$ are
homing into the generating objects $\alg{2}$ and $\utilde{\alg{2}}$,
respectively.  This situation is an instance of a more general framework
investigated at the end of the 20th century in the universal algebra community
under the name of natural dualities.

\subsection{Natural dualities in a nutshell}\label{subsec:nutshell}

The general setting of natural dualities involves a finite algebra $\alg{N}$ and
a \emph{topological alter ego} $\utilde{\alg{N}}$, which carries the same
underlying set $N$ equipped with the discrete topology and possibly additional
relations and (partial) operations. The key requirement is that this extra
structure be \emph{compatible} with $\alg{N}$, meaning that every relation and
every (partial) operation of $\utilde{\alg{N}}$ is recognizable from the
algebraic side. Specifically, every $n$-ary relation must be a subalgebra of $\alg{N}^n$,
every $n$-ary total operation must be a homomorphism $\alg{N}^n \to \alg{N}$,
and every $n$-ary partial operation must have a subalgebra of $\alg{N}^n$ as
its domain and be a homomorphism from that subalgebra to $\alg{N}$.

The quasivariety $\mathcal{A} = \classop{ISP}(\alg{N})$ and the category
$\mathcal{X} = \classop{IS}_c\classop{P}^+(\utilde{\alg{N}})$ of isomorphic
copies of closed substructures of powers of $\utilde{\alg{N}}$ over nonempty
index sets are then related by two
hom-functors. For each $\alg{A} \in \mathcal{A}$, the hom-set
$\mathcal{A}(\alg{A}, \alg{N})$ inherits a structure from the product
$\utilde{\alg{N}}^A$ via pointwise operations and the subspace topology, yielding a functor
\[
  D \colon \mathcal{A} \longrightarrow \mathcal{X}, \quad \alg{A} \mapsto \mathcal{A}(\alg{A}, \alg{N}).
\]
On morphisms, the functor $D$ acts by precomposition, sending a homomorphism
$f\colon \alg{A} \to \alg{B}$
to $D(f)\colon D(\alg{B}) \to D(\alg{A})$, $h \mapsto h \circ f$.
Symmetrically, for each $\alg{X} \in \mathcal{X}$, the hom-set
$\mathcal{X}(\alg{X}, \utilde{\alg{N}})$ of continuous structure-preserving maps
inherits an algebraic structure from $\alg{N}^X$ via pointwise operations, yielding a functor
\[
  E \colon \mathcal{X} \longrightarrow \mathcal{A}, \quad \alg{X} \mapsto \mathcal{X}(\alg{X}, \utilde{\alg{N}}).
\]
On morphisms, the functor $E$ acts by precomposition, sending a continuous
structure-preserving map
$g\colon \alg{X} \to \alg{Y}$ to $E(g)\colon E(\alg{Y}) \to E(\alg{X})$,
$k \mapsto k \circ g$.
The pair $(D, E)$ forms a \emph{dual adjunction} between $\mathcal{A}$ and
$\mathcal{X}$, with natural transformations
\[
  e_{\alg{A}} \colon \alg{A} \to ED(\alg{A}), \qquad
  \varepsilon_{\alg{X}} \colon \alg{X} \to DE(\alg{X}),
\]
both given by evaluation, $e_{\alg{A}}(a)(h) = h(a)$ and
$\varepsilon_{\alg{X}}(x)(k) = k(x)$. Following \cite[\S2]{ClarkDavey} we call
them the \emph{evaluation maps}, and avoid the words unit and counit.

Both evaluation maps are always embeddings, and the theory is concerned with
when they are isomorphisms. One says that
$\utilde{\alg{N}}$ \emph{yields a duality} on $\mathcal{A}$ if $e_{\alg{A}}$ is
an isomorphism for every $\alg{A} \in \mathcal{A}$. The duality is
\emph{full} if moreover $\varepsilon_{\alg{X}}$ is an isomorphism for every
$\alg{X} \in \mathcal{X}$, in which case $D$ and $E$ are
dual equivalences. Requiring these two conditions for the finite objects of
$\mathcal{A}$ and of $\mathcal{X}$ only defines a \emph{finite-level duality}
and a \emph{finite-level full duality}, respectively. Finally, the duality is
\emph{strong} if it is full and, in addition, the alter ego $\utilde{\alg{N}}$
is injective in $\mathcal{X}$, meaning that whenever $\alg{Y}$ is a closed
substructure of an object $\alg{X}$ of $\mathcal{X}$, every morphism $\alg{Y}
\to \utilde{\alg{N}}$ extends to a morphism $\alg{X} \to \utilde{\alg{N}}$.
Two consequences are used below without mention. Precomposition with a
surjective homomorphism is an embedding \cite[Theorem~1.5.3(iii)]{ClarkDavey},
and, when the duality is strong, a morphism of $\mathcal{X}$ is an embedding if
and only if precomposition with it is surjective
\cite[Lemma~3.2.6]{ClarkDavey}.
The theory of natural dualities investigates, among other things, sufficient
conditions on $\utilde{\alg{N}}$ for the dual adjunction to be a duality, a full
duality, or a strong duality. Its two motivating examples are Stone duality for
Boolean algebras, where $\utilde{\alg{N}}$ is the discrete two element space
\cite[Theorem~4.1.2]{ClarkDavey}, and Priestley duality for distributive
lattices, where $\utilde{\alg{N}}$ is the discrete two element chain with its
bounds as constants \cite[Theorem~1.2.5]{ClarkDavey}.

In our setting, $\alg{N} = \alg{2}$ is the two element median algebra and
$\utilde{\alg{N}} = \utilde{\alg{2}}_\tau$ is the two element bounded poc set
equipped with the discrete topology. We denote the functors $D$ and $E$ by $-^*$
and $-_*$, respectively. The source category $\mathcal{A}$ specializes to
$\cate{M} = \classop{ISP}(\alg{2})$, and the dual category $\mathcal{X}$
specializes to the category $\Poc_\tau :=
\classop{IS}_{c}\classop{P}^+(\utilde{\alg{2}}_\tau)$ of \emph{Boolean
topological bounded poc sets}, whose objects are the topologically
closed sub-poc sets of powers of
$\utilde{\alg{2}}_\tau$, equipped with the product topology. Here
$\classop{P}^+$ denotes products over nonempty index sets. A subscript $0$
on $\classop{S}$, which will reappear in the TopSwap Theorem below, allows the
empty substructure. This refinement is vacuous in our
setting,\footnote{No substructure is empty, because the type of
$\utilde{\alg{2}}$ contains the constants $0$ and $1$.} whereas the
restriction to nonempty index sets excludes the empty power, a one element
structure in which $0$ and $1$ coincide. Forgetting the topology, we may
therefore write
$\classop{IS}_0\classop{P}^+(\utilde{\alg{2}}) =
\classop{ISP}^+(\utilde{\alg{2}}) = \cate{Poc}$ by
Proposition~\ref{prop:char}.

The objects of $\Poc_\tau$ are also called \emph{bounded strongly complemented
Priestley spaces} \cite{ClarkDavey}. They are the bounded Priestley
spaces carrying an order reversing involutive homeomorphism $(-)^\co$ that
swaps the bounds and satisfies $x \leq x^\co$ only if $x = 0$.

\subsection{Two dualities from one machine}\label{subsec:two-dualities}

Everything we need from the theory is contained in two of its results, both
stated below in the notation of Subsection~\ref{subsec:nutshell}. The first one
produces a strong duality out of a single syntactic hypothesis on $\alg{N}$,
the existence of a near-unanimity term operation, and the second one turns a
duality already obtained into its partner, with the topology moved to the other
side.

\begin{definition}\label{def:nu}
An $n$-ary operation $t$ on a set $N$ is called a \emph{near-unanimity} (NU)
\emph{operation} if
\[
  t(x, \ldots, x, \underset{i}{y}, x, \ldots, x) = x
\]
holds for every position $i$ of $y$. A ternary NU operation is called a
\emph{majority operation}.
\end{definition}

\begin{theorem}\label{thm:nd-toolkit}
Let $\alg{N}$ be a finite algebra of finite type and let $\utilde{\alg{N}}$ be
an alter ego of $\alg{N}$.
\begin{enumerate}
\item\label{it:nu} \textup{(NU Strong Duality Corollary
{\cite[Corollary~3.3.9]{ClarkDavey}})} If $\alg{N}$ has a majority term
operation and
every non-trivial subalgebra of $\alg{N}$ is subdirectly irreducible, then the
alter ego
\[
  \utilde{\alg{N}}_0 := \langle N,\; K,\; P_1,\; \classop{S}(\alg{N}^2),\;
  \tau_d \rangle
\]
yields a strong duality on $\mathcal{A}$. Here $K$ is the union of
the one element subalgebras of $\alg{N}$, the set $P_1$ consists of the unary
partial endomorphisms of $\alg{N}$, and $\classop{S}(\alg{N}^2)$ is the set
of subalgebras of $\alg{N}^2$, read as binary relations.
\item\label{it:topswap} \textup{(TopSwap Theorem
{\cite[Theorem~2.4]{DaveyHaviarPriestley2012}})} If $\utilde{\alg{N}}$ yields a
finite-level full duality between $\mathcal{A}$ and $\mathcal{X}$, the empty
substructure being allowed on the dual side, then the associated hom-functors
yield a dual equivalence between $\classop{IS}_c\classop{P}(\alg{N}, \tau_d)$
and $\classop{IS}_0\classop{P}^+(\utilde{\alg{N}})$, where the subscript $0$
allows the empty substructure, the first category carrying
the topology on the algebraic side and the second one forgetting it.
\end{enumerate}
\end{theorem}

Both hypotheses of item~(\ref{it:nu}) are met by $\alg{2}$. The median
operation $m$ is a majority operation, since $m(y,x,x) = m(x,y,x) = m(x,x,y) =
x$ by definition, and $\alg{2}$ has no proper non-trivial subalgebra. The alter
ego that item~(\ref{it:nu}) supplies is larger than the one we wish to work
with, and the next lemma trims it. The tool for such trimming is the notion of
\emph{entailment} of \cite[\S2.4]{ClarkDavey}, under which an alter ego entails
a second one on the same underlying set when every map preserving the structure
of the first preserves that of the second as well. We use the fact that a map
preserves
a (partial) operation if and only if it preserves the graph of that operation,
read as a relation, which lets us argue
with relations alone.

\begin{lemma}\label{lem:same-morphisms}
The alter egos $\utilde{\alg{2}}_\tau$ and $\utilde{\alg{2}}_0$ entail each
other in the sense of \cite[2.4.1]{ClarkDavey}. They determine the same dual
category, the same pair of hom-functors and the same evaluation maps.
\end{lemma}
\begin{proof}
By \cite[Lemma~2.4.4]{ClarkDavey} it is enough, in each direction, to check that
every morphism from a closed substructure of a power into the one alter ego
preserves the structure of the other.

Each item of structure of $\utilde{\alg{2}}_\tau$ occurs in
$\utilde{\alg{2}}_0$. The bounds $0$ and $1$ are the one element subalgebras of
$\alg{2}$, the complementation $(-)^\co$ is an endomorphism of $\alg{2}$ and
hence a member of $P_1$, and the order $\leq$ is a subalgebra of $\alg{2}^2$.
Every morphism for $\utilde{\alg{2}}_0$ is therefore a morphism for
$\utilde{\alg{2}}_\tau$.

For the converse, observe first that every nonempty subset of $\{0,1\}^2$ is a
subalgebra of $\alg{2}^2$. There are thus fifteen binary algebraic relations on
$\{0,1\}$, and each of them is definable from $\leq$, $(-)^\co$, $0$ and $1$ by
a conjunction of atomic formulas, hence preserved by every morphism for
$\utilde{\alg{2}}_\tau$. The whole square is defined by $x = x$, the relations
$\leq$ and $\geq$ by $x \leq y$ and by $y \leq x$, the two remaining three
element subsets by $x \leq y^\co$ and by $x^\co \leq y$, the diagonal by $x
\leq y \wedge y \leq x$, the antidiagonal by $y = x^\co$, the four two element
subsets that fix a coordinate by $x = 0$, $x = 1$, $y = 0$ and $y = 1$, and
each singleton by fixing both coordinates. Since $K = \{0,1\}$ and every member
of $P_1$ has its graph among the fifteen relations just listed, every morphism
for $\utilde{\alg{2}}_\tau$ is a morphism for $\utilde{\alg{2}}_0$.

Mutual entailment settles the morphisms. The two alter egos also have the
same closed substructures in every power, since $P_1$ consists of the
identity, the complementation, the two constant maps and their restrictions
to the one element subalgebras, while relations impose no closure condition.
\end{proof}

The second duality moves the topology to the algebraic side. Its source
category is
\[
  \cate{M}_\tau := \classop{IS}_c\classop{P}(\alg{2}_\tau),
\]
the category of \emph{Boolean topological median algebras}, whose objects are
topologically closed subalgebras of powers of $\alg{2}_\tau$ equipped with the
product topology, and whose morphisms are continuous median algebra homomorphisms.
Its two hom-functors are defined on objects by
\begin{align*}
  F \colon \cate{M}_\tau \longrightarrow \cate{Poc},&\quad
    F(\alg{A}) = \cate{M}_\tau(\alg{A}, \alg{2}_\tau),\\
  {-}_+ \colon \cate{Poc} \longrightarrow \cate{M}_\tau,&\quad
    \alg{X}_+ = \cate{Poc}(\alg{X}, \utilde{\alg{2}}),
\end{align*}
where $F(\alg{A})$ carries the structure of a substructure of
$\utilde{\alg{2}}^A$ (making it a bounded poc set) and $\alg{X}_+$ carries the
closed subalgebra structure of $\alg{2}_\tau^X$ (making it a Boolean topological
median algebra). Both act on morphisms by precomposition. The functor $-_+$
here is the complex algebra construction, now equipped with the Boolean topology
inherited from $\alg{2}_\tau^X$. By contrast, the set $F(\alg{A})$ consists only
of
the \emph{continuous} homomorphisms $\alg{A} \to \alg{2}$ (in general a
strict subset of $\alg{A}^+$).

\begin{theorem}\label{thm:dualities}
\begin{enumerate}
\item\label{it:dual-poctau} The functors $-^*$ and $-_*$ yield a strong duality
between $\cate{M}$ and $\Poc_\tau$.
\item\label{it:dual-poc} The functors $F$ and $-_+$ yield a dual equivalence
between $\cate{M}_\tau$ and $\cate{Poc}$.
\end{enumerate}
\end{theorem}
\begin{proof}
For~(\ref{it:dual-poctau}), recall that $\alg{2}$ satisfies both
assumptions of Theorem~\ref{thm:nd-toolkit}(\ref{it:nu}), so that
$\utilde{\alg{2}}_0$ yields a strong duality on $\cate{M} =
\classop{ISP}(\alg{2})$. The statement now follows from
Lemma~\ref{lem:same-morphisms}.

For~(\ref{it:dual-poc}), a strong duality is full, and a full duality is in
particular full at the finite level, so~(\ref{it:dual-poctau}) supplies
Theorem~\ref{thm:nd-toolkit}(\ref{it:topswap}) with its hypothesis, and the
latter provides the dual equivalence of the statement.
\end{proof}

Both dualities are known in geometric group theory. Roller states them as
anti-equivalences of categories \cite[Theorems~5.3 and~6.4]{Roller2016} and
defines the dual maps \cite[\S7]{Roller2016}, while Bowditch
\cite[\S9.2, Theorem~12.5.1 and \S26.12]{Bowditch2022} works at the level of
objects.

The two dualities are not independent outputs of the machine. One
turn of the crank produces the first one, and the theory then hands over the
second one for free, asking only for the fullness of the first.

\subsection{Boolean topological median algebras as profinite algebras}
The objects of $\cate{M}_\tau$ admit a characterization as profinite
algebras. Neither the notion nor the constructions behind it are specific to
median algebras. Profinite algebras and profinite completions form a standard
chapter of universal algebra \cite{Numakura1957, Banaschewski1972,
ClarkDaveyFreeseJackson2004,
DaveyGouveiaHaviarPriestley2011, HansoulTeheux}. We recall the relevant
constructions, applying them to the variety $\cate{M}$. A \emph{co-directed diagram} of median algebras is a family
$((\alg{B}_i)_{i \in I},\, (\varphi_{ij})_{i \leq j})$ indexed by a
co-directed poset $I$ (a poset in which any two elements have a lower bound),
where each bonding map $\varphi_{ij} \colon \alg{B}_i \to \alg{B}_j$ is a
median algebra homomorphism, and the conditions $\varphi_{ii} =
\mathrm{id}_{\alg{B}_i}$ and $\varphi_{jk} \circ \varphi_{ij} = \varphi_{ik}$
hold for all $i \leq j \leq k$ in $I$.

A \emph{cone} to such a diagram is a median algebra $\alg{C}$ together with
homomorphisms $f_i \colon \alg{C} \to \alg{B}_i$ satisfying $\varphi_{ij} \circ
f_i = f_j$ for all $i \leq j$ in $I$. The \emph{limit} $\alg{B}:=\varprojlim
\{\alg{B}_i\mid i \in I\}$  is the terminal cone, equipped with cone maps $p_i
\colon \alg{B} \to \alg{B}_i$. By definition, for every cone $(\alg{C}, (f_i)_{i \in I})$ to
the original diagram, there is a unique homomorphism $h \colon \alg{C} \to
\alg{B}$ satisfying $p_i \circ h = f_i$ for all $i \in I$. When each $\alg{B}_i$
is a topological median algebra, the same construction applies with continuous homomorphisms
throughout, and the limit carries the initial topology with respect to the cone
maps $(p_i)_{i \in I}$.

The limit always exists and it is the subalgebra
\[
  \varprojlim \{\alg{B}_i \mid i \in I\}
  \;=\;
  \bigl\{\, b \in \textstyle\prod_{i \in I} \alg{B}_i
  \;\bigm|\;
  \varphi_{ij}(b_i) = b_j \text{ for all } i \leq j \,\bigr\}
\]
of the direct product of the $\alg{B}_i$, the restriction of the projection maps
serving as cone maps. Its elements are called \emph{compatible threads}.  When
each $\alg{B}_i$ is a Hausdorff topological median algebra, each condition
$\varphi_{ij}(b_i) = b_j$ is closed (it is the equalizer of two continuous maps
into a Hausdorff space), so the limit is a closed subalgebra of the product
$\prod_{i \in I} \alg{B}_i$. In particular, when each $\alg{B}_i$ is finite and
discrete, the product is compact, so the limit is compact,
Hausdorff, and totally disconnected.

\begin{definition}\label{def:profinite}
A topological median algebra $\alg{B}$ is \emph{profinite} if it is isomorphic
algebraically and topologically to the limit of a co-directed diagram of finite
discrete median algebras.
\end{definition}

In a product of discrete spaces a basic open set constrains finitely many
coordinates, so an element $c$ of $\prod_{i \in I} X_i$ lies in the closure of a
subset $B$ if and only if every finite $F \subseteq I$ (in notation $F \Subset
I$) admits an element $b \in B$ with $\pi_F(b) = \pi_F(c)$. Following tradition,
we then say that $c$ is \emph{locally in} $B$, and we call such an element $b$ a
\emph{witness} for $c$ on $F$.  When $B$ is the image of an embedding $\iota
\colon \alg{A} \hookrightarrow \prod_{i \in I} X_i$, we transfer the vocabulary
to the domain, saying that $c$ is locally in $\alg{A}$ and that $a \in A$
witnesses $c$ on $F$. Ordered by reverse inclusion, the finite subsets of $I$
index a co-directed diagram of finite images $\pi_F(B)$ whose bonding maps are
the coordinate projections, and one can check that
\begin{equation}\label{eq:locality}
  \overline{B} \;\cong\;
  \varprojlim \{\pi_F(B) \mid F \Subset I\}.
\end{equation}
When the $X_i$ are median algebras and $B$ is a subalgebra, the two objects
are median algebras and the \emph{locality identification}~\eqref{eq:locality}
is an isomorphism in $\cate{M}_\tau$. In particular a closed subalgebra of such
a product is the limit of its finite images.

Our characterization of $\cate{M}_\tau$ relies on well-established tools from
universal algebra. A congruence $\theta$ of an algebra $\alg{A}$
\emph{saturates} a subset $U$ of $A$ when $U$ is a union of $\theta$-classes,
and the \emph{syntactic congruence} of $U$ is the largest congruence saturating
$U$. It always exists, and it is computed by declaring $a$ and $b$ equivalent
when $p(a) \in U$ if and only if $p(b) \in U$, for every term operation $p$ in
one variable and with parameters from $\alg{A}$. A variety has \emph{finitely
determined syntactic congruences}, abbreviated to FDSC, when a fixed finite set
of terms is enough in that computation, uniformly over the members of the
variety and over the subsets \cite[\S2]{ClarkDaveyFreeseJackson2004}. Median
algebras satisfy this in the strongest possible way, one term with two
parameters being enough.

\begin{lemma}\label{lem:median-fdsc}
The variety of median algebras has FDSC. For any subset $U$ of a median algebra
$\alg{A}$, the syntactic congruence of $U$ is the relation $\sim_U$ defined by
\[
  a \sim_U b \quad\text{whenever}\quad
  \forall u, v \in A \quad (m(a,u,v) \in U \iff  m(b,u,v) \in U).
\]
\end{lemma}
\begin{proof}
Write $T(a) := \{(u,v) \in A^2 \mid m(a,u,v) \in U\}$, so that $a \sim_U b$
means $T(a) = T(b)$. The relation $\sim_U$ is a congruence, because the
distributivity axiom \eqref{eq:distributivity} gives
\[
  m\bigl(m(a,r,s),u,v\bigr) = m\bigl(a, m(r,u,v), m(s,u,v)\bigr)
\]
for all $r$, $s$, $u$ and $v$ in $A$, and the right-hand side depends on $a$
only through $T(a)$, so that $a \sim_U b$ forces $m(a,r,s) \sim_U m(b,r,s)$.

The relation $\sim_U$ saturates $U$. Indeed, let $a \in U$ and $b \notin U$. The
pair $(a,b)$ belongs to $T(a)$ but not to $T(b)$, since $m(a,a,b) = a$ lies in
$U$ whereas $m(b,a,b) = b$ does not.

Finally, let $\theta$ be a congruence saturating $U$ and let $a \mathrel{\theta}
b$. For all $u$ and $v$ in $A$ the elements $m(a,u,v)$ and $m(b,u,v)$ are
$\theta$-equivalent, hence they lie both in $U$ or both outside it, so $T(a) =
T(b)$. 
\end{proof}

\begin{proposition}\label{prop:profinite-char}
For a topological median algebra $\alg{B}$, the following conditions are equivalent.
\begin{enumerate}
\item\label{it:prof01} The algebra $\alg{B}$ is profinite.
\item\label{it:prof02} The algebra $\alg{B}$ is Boolean topological.
\item\label{it:prof03} The algebra $\alg{B}$ belongs to $\cate{M}_\tau$.
\end{enumerate}
Moreover, products of profinite median algebras are profinite, and closed
subalgebras of profinite median algebras are profinite.
\end{proposition}
\begin{proof}
(\ref{it:prof03})$\Rightarrow$(\ref{it:prof02}).
Every object of $\cate{M}_\tau$ is a closed subalgebra of a power of
$\alg{2}_\tau$, hence Boolean topological.

(\ref{it:prof02})$\Rightarrow$(\ref{it:prof01}).
The variety of median algebras has FDSC by Lemma~\ref{lem:median-fdsc}, and a
Boolean topological algebra whose underlying algebra lies in a variety with FDSC
is profinite \cite[Theorem~8.1]{ClarkDaveyFreeseJackson2004}.

(\ref{it:prof01})$\Rightarrow$(\ref{it:prof03}).  Let $\alg{B}$ be the limit of
a co-directed diagram of finite median algebras $\alg{B}_i$ with the discrete
topology. By
Theorem~\ref{prop:med-axioms} each $\alg{B}_i$ embeds into the power
$\alg{2}_\tau^{\M(\alg{B}_i, \alg{2})}$. Hence $\prod_i \alg{B}_i$ is a closed
subalgebra of a power of $\alg{2}_\tau$, and so is its closed subalgebra
$\alg{B}$. Thus $\alg{B} \in \cate{M}_\tau$.

The stability statement follows from the equivalence with
condition~(\ref{it:prof03}), since $\cate{M}_\tau =
\classop{IS}_c\classop{P}(\alg{2}_\tau)$ is closed under $\classop{P}$ and
$\classop{S}_c$.
\end{proof}

By Proposition~\ref{prop:profinite-char}, every Boolean topological median
algebra embeds in a power of $\alg{2}_\tau$. Using the coordinate projections
to separate points, the following separation property comes for free.

\begin{corollary}[Clopen separation {\cite[Proposition~5.4]{Roller2016}}]\label{lem:clopen-sep}
Let $\alg{B}$ be a Boolean topological median algebra. For any two distinct
elements $x$ and $y$ of $\alg{B}$, there is a clopen prime convex subset of
$\alg{B}$ that contains $x$ but not $y$.
\end{corollary}

\begin{remark}\label{rem:standard}
Median algebras also inherit FDSC from a general fact, since every majority term
makes a variety congruence-distributive, and every finitely generated
congruence-distributive variety has FDSC
\cite[\S6]{ClarkDaveyFreeseJackson2004}, see also
\cite{DaveyJacksonMarotiMcKenzie2008}. Lemma~\ref{lem:median-fdsc} is more
precise, and it explains Roller's proof of Corollary~\ref{lem:clopen-sep},
which uses the syntactic congruence of a clopen set.
\end{remark}

Combining Proposition~\ref{prop:profinite-char} with
Theorem~\ref{thm:dualities}(\ref{it:dual-poc}) records the duality in the form
we shall use.

\begin{corollary}\label{cor:profinite-duality}
The functors $F$ and $-_+$ yield a dual equivalence between the category of
profinite median algebras with continuous homomorphisms and the category
$\cate{Poc}$ of bounded poc sets.
\end{corollary}

\section{Profinite completions of median
algebras}\label{sec:profinite}

We recall a standard topological fact used repeatedly in what follows.
If $f, g \colon X \to Y$ are continuous maps and $Y$ is Hausdorff, then
$f = g$ whenever the two maps agree on a dense subset of $X$.

Denote by $-^\flat \colon \Poc_\tau \to \Poc$ the forgetful functor that
discards the topology. The composition $-_+ \circ -^\flat \circ -^*$ interleaves
the two dualities, and this composition reduces to $-_+ \circ -^+$. As
Corollary~\ref{prop:profinite} shows, the resulting object $(\alg{A}^+)_+$ is
the profinite completion of $\alg{A}$.

\subsection{Profinite completion}

The profinite completion is a classical construction of universal algebra
\cite{Numakura1957,Banaschewski1972} that has become a tool of algebraic logic,
where it computes canonical extensions and universal models
\cite{Harding2006,Bezhanishvili2006,Vosmaer2010,Gehrke2014,%
DeBerardinisGhilardi2024}.  We recall the construction and instantiate it for
median algebras.

Let $\mathcal{S}_{\alg{A}}$ denote the family of all congruences of finite
index on $\alg{A}$, ordered by inclusion. Any two of its
members $\alpha$ and $\beta$ have the lower bound $\alpha \cap \beta$, again of
finite index, so the family is co-directed. For $\alpha \subseteq \beta$, the natural surjection
$\varphi_{\alpha\beta} \colon \alg{A}/\alpha \to \alg{A}/\beta$ serves as bonding
map, and $(\alg{A}/\alpha, \varphi_{\alpha\beta})_{\alpha \in
\mathcal{S}_{\alg{A}}}$ becomes a co-directed diagram of finite median
algebras with the discrete topology.

Define the \emph{profinite completion} of $\alg{A}$ as the limit of this diagram,
\[
  \hat{\alg{A}} := \varprojlim \{\alg{A}/\alpha \mid \alpha \in
  \mathcal{S}_{\alg{A}}\}.
\]
By construction, the completion $\hat{\alg{A}}$ is profinite, hence an object
of $\cate{M}_\tau$ by Proposition~\ref{prop:profinite-char}. The quotient
maps $\pi_\alpha\colon \alg{A} \to \alg{A}/\alpha$ form a cone over the diagram and induce the
canonical homomorphism
\[
  \mu_{\alg{A}} \colon \alg{A} \longrightarrow \hat{\alg{A}}, \qquad a \mapsto
  (a/\alpha)_{\alpha \in \mathcal{S}_{\alg{A}}}.
\]
We write $p_\alpha \colon \hat{\alg{A}} \to \alg{A}/\alpha$ for the limit
projections, so that $p_\alpha \circ \mu_{\alg{A}} = \pi_\alpha$.

The homomorphism $\mu_{\alg{A}}$ is injective. Indeed, for distinct elements
$a$ and $b$ of $A$, there is a homomorphism $h \in \alg{A}^*$ with
$h(a) \neq h(b)$, and the
kernel of $h$ belongs to $\mathcal{S}_{\alg{A}}$ and separates their images.

The image $\mu_{\alg{A}}(A)$ is dense in $\hat{\alg{A}}$, as any basic open set
of $\hat{\alg{A}}$ fixes finitely many coordinates $\alpha_1, \ldots,
\alpha_n$, and their meet $\alpha = \alpha_1 \cap \cdots \cap \alpha_n$ again lies in
$\mathcal{S}_{\alg{A}}$. Any element $a \in A$ representing the
$\alpha$-coordinate of a given thread agrees with that thread at each $\alpha_i$
by compatibility of the bonding maps. Hence $\mu_{\alg{A}}(a)$ meets the chosen
open set. Consequently $\hat{\alg{A}}$ is the closure of the image,
\[
  \hat{\alg{A}} = \overline{\mu_{\alg{A}}(A)} \;\subseteq\; \prod_{\alpha \in
  \mathcal{S}_{\alg{A}}} \alg{A}/\alpha,
\]
and the locality identification~\eqref{eq:locality} identifies the compatible
threads of the diagram with the elements of the product that are locally in
$\alg{A}$ with respect to $\mu_{\alg{A}}$.

Let $f \colon \alg{A} \to \alg{B}$ be a homomorphism in $\cate{M}$. For each
$\beta \in \mathcal{S}_{\alg{B}}$, the pullback $f^{-1}(\beta)$ is a congruence of
finite index on $\alg{A}$, and $f$ induces an embedding $\iota_\beta \colon
\alg{A}/f^{-1}(\beta) \hookrightarrow \alg{B}/\beta$. Composing the projection
$p_{f^{-1}(\beta)}$ with this embedding yields continuous homomorphisms
$\hat{\alg{A}} \to \alg{B}/\beta$ compatible with the bonding maps of $\alg{B}$,
so that $\hat{\alg{A}}$ becomes a cone over the diagram defining $\hat{\alg{B}}$.
Since $\hat{\alg{B}}$ is the terminal cone over that diagram, this cone factors
through a unique continuous homomorphism $\hat{f} \colon \hat{\alg{A}} \to
\hat{\alg{B}}$, namely $x \mapsto
\bigl(\iota_\beta(x_{f^{-1}(\beta)})\bigr)_{\beta \in \mathcal{S}_{\alg{B}}}$,
making $\hat{-} \colon \cate{M} \to \cate{M}_\tau$ a functor.

The limit cone of the preceding construction characterizes $\hat{\alg{A}}$ by
its maps into the finite quotients $\alg{A}/\alpha$.
The next proposition points the other way, describing the maps out of $\alg{A}$
into an arbitrary Boolean topological target. Write
\[
  U \colon \cate{M}_\tau \longrightarrow \cate{M}
\]
for the functor that forgets the topology, which is well defined because a
closed subalgebra of $\alg{2}_\tau^I$ is in particular a subalgebra of
$\alg{2}^I$.

The proposition that follows is not specific to median algebras. For any
variety, the profinite completion is left adjoint to
the functor forgetting the topology \cite[Proposition~8 and
Corollary~1]{Banaschewski1972}. The proof below is given because it is explicit
about the extension.

\begin{proposition}[The profinite adjunction]
\label{prop:ump-reflection}
The profinite completion functor $\hat{-}$ is left adjoint to $U$. In other
words, for every $\alg{A} \in \cate{M}$ and $\alg{B} \in \cate{M}_\tau$ there
is a bijection
\[
  \cate{M}\bigl(\alg{A}, U(\alg{B})\bigr) \;\cong\;
  \cate{M}_\tau\bigl(\hat{\alg{A}}, \alg{B}\bigr),
\]
natural in $\alg{A}$ and $\alg{B}$. The unit of the adjunction is
$\mu_{\alg{A}}$.
\end{proposition}
\begin{proof}
Every continuous homomorphism $h \colon \hat{\alg{A}} \to \alg{B}$ restricts
along $\mu_{\alg{A}}$ to the homomorphism $U(h) \circ \mu_{\alg{A}} \colon
\alg{A} \to U(\alg{B})$, and the bijection of the statement is this restriction
map. We have to show that every homomorphism $f \colon \alg{A} \to U(\alg{B})$
extends to exactly one continuous homomorphism $f^\sharp \colon \hat{\alg{A}}
\to \alg{B}$, in the sense that $f = U(f^\sharp) \circ \mu_{\alg{A}}$.
Uniqueness is the injectivity of the restriction map, and existence its
surjectivity.

\textit{Uniqueness.} Any two continuous extensions of $f$ agree on
$\mu_{\alg{A}}(A)$, which is dense in $\hat{\alg{A}}$. Since $\alg{B}$ is
Hausdorff, uniqueness follows.

\textit{Existence.} Realize $\alg{B}$ as a closed subalgebra of
$\alg{2}_\tau^I$, and write $\rho_i \colon \alg{2}_\tau^I \to \alg{2}_\tau$ for
the coordinate projections. For $i \in I$, the composite $f_i := \rho_i \circ f
\colon \alg{A} \to \alg{2}$ is a homomorphism into the two element algebra, so
its kernel $\alpha_i$ lies in
$\mathcal{S}_{\alg{A}}$, and $f_i$ factors as $u_i \circ \pi_{\alpha_i}$ for a
homomorphism $u_i \colon \alg{A}/\alpha_i \to \alg{2}$. The quotient
$\alg{A}/\alpha_i$ is finite and discrete, so $u_i$ is continuous, and the limit
projection $p_{\alpha_i} \colon \hat{\alg{A}} \to \alg{A}/\alpha_i$ is
continuous as one of the maps defining the topology of $\hat{\alg{A}}$. Each
composite $u_i \circ p_{\alpha_i}$ is therefore
continuous, and together these maps define a homomorphism $g \colon
\hat{\alg{A}} \to \alg{2}_\tau^I$, continuous because the product topology is
initial with respect to the family $(\rho_i)_{i \in I}$. For $a \in A$ and $i
\in I$ we have
\[
  \rho_i\bigl(g(\mu_{\alg{A}}(a))\bigr) = u_i\bigl(p_{\alpha_i}(\mu_{\alg{A}}(a))\bigr)
  = u_i(a/\alpha_i) = f_i(a),
\]
that is, the identity $g \circ \mu_{\alg{A}} = f$ holds. In particular $g$ maps
$\mu_{\alg{A}}(A)$ into $\alg{B}$, so the closed set $g^{-1}(\alg{B})$ contains
$\mu_{\alg{A}}(A)$, which is dense, and therefore contains all of
$\hat{\alg{A}}$. Corestricting $g$ to $\alg{B}$ gives $f^\sharp$.

\textit{Naturality.} The naturality of the bijection in $\alg{B}$ is the
functoriality of $U$. Its naturality in $\alg{A}$, along $t \colon \alg{A}'
\to \alg{A}$, amounts to the identity
\begin{equation}\label{eq:mu-natural}
  \hat{t} \circ \mu_{\alg{A}'} = \mu_{\alg{A}} \circ t ,
\end{equation}
which states that the canonical homomorphisms are natural, or equivalently that
$\hat{t}$ extends $t$. To verify it, recall that $\hat{t}$ is determined by the
identities $p_\alpha \circ \hat{t} = \iota_\alpha \circ p_{t^{-1}(\alpha)}$,
for $\alpha \in \mathcal{S}_{\alg{A}}$, where $\iota_\alpha \colon
\alg{A}'/t^{-1}(\alpha) \hookrightarrow \alg{A}/\alpha$ is the embedding
induced by $t$. Evaluating them at $\mu_{\alg{A}'}(a)$, for $a \in A'$, gives
\[
  p_\alpha\bigl(\hat{t}(\mu_{\alg{A}'}(a))\bigr)
  = \iota_\alpha\bigl(a/t^{-1}(\alpha)\bigr)
  = t(a)/\alpha
  = p_\alpha\bigl(\mu_{\alg{A}}(t(a))\bigr)
\]
for every $\alpha$, and the projections are jointly injective.
\end{proof}

In practice we use the proposition in the following form. Any median algebra
homomorphism $f \colon \alg{A} \to \alg{B}$ from a median algebra $\alg{A}$
into a profinite one admits a unique continuous extension $f^\sharp \colon
\hat{\alg{A}} \to \alg{B}$ satisfying $f^\sharp \circ \mu_{\alg{A}} = f$.

When $\alg{B}$ already carries a topology, the family $\mathcal{S}_{\alg{B}}$
still consists of all congruences of finite index on the abstract algebra
$U(\alg{B})$, and such a congruence need not be open. The canonical
homomorphism $\mu_{\alg{B}}$ may then fail even to be continuous, as the
following example shows.

\begin{example}\label{ex:not-strongly-complete}
Let $\alg{B} = \alg{2}_\tau^{\N}$, a profinite median algebra whose elements
we identify with the subsets of $\N$, so that
$m(X,Y,Z) = (X \cap Y) \cup (X \cap Z) \cup (Y \cap Z)$. Given an ultrafilter
$\mathcal{V}$ on $\N$, let $u \colon \alg{B} \to \alg{2}$ send $Z$ to $1$ if
and only if $Z \in \mathcal{V}$. It is a lattice homomorphism from the power
set of $\N$ onto $\{0,1\}$, hence a homomorphism of median algebras, the
median being a lattice term. Its kernel $\beta$ is therefore a congruence of
index two whose classes are $\mathcal{V}$ and $\mathcal{V}^\co$, and
$\mathcal{V}$ is a prime convex subset of $\alg{B}$, see also
\cite[Example~2.24']{Fioravanti2020}.

Assume that $\mathcal{V}$ is non-principal. Then $u$ is not continuous.
Otherwise $\mathcal{V}$ would be open and, being nonempty, would contain a
basic clopen set $\{Z \subseteq \N \mid Z \cap F = S\}$, with $F \Subset \N$
and $S \subseteq F$. That set contains the finite set $S$, whereas a
non-principal ultrafilter contains no finite set. Up to the isomorphism
$\alg{B}/\beta \cong \alg{2}$, the quotient map $\pi_\beta$ is $u$, a
homomorphism of median algebras that is not continuous, so the functor $U$ is
not full. Since $p_\beta$ is continuous and $p_\beta \circ \mu_{\alg{B}} =
\pi_\beta$, the map $\mu_{\alg{B}}$ is not continuous either.
\end{example}

The next corollary identifies the median algebras for which this does not
happen.

\begin{corollary}\label{cor:compact-self-complete}
Let $\alg{A}$ be a median algebra, and endow $A$ with the topology induced by
$\mu_{\alg{A}}$ from $\hat{\alg{A}}$. Then $\mu_{\alg{A}}$ is onto if and only
if $\alg{A}$ is compact, in which case it is an isomorphism of
$\cate{M}_\tau$.
\end{corollary}
\begin{proof}
The forward implication is clear since $\hat{\alg{A}}$ is compact. Conversely,
assume that $\alg{A}$ is
compact. Its image is then a compact subspace of the Hausdorff space
$\hat{\alg{A}}$, hence closed, and it is dense, hence the whole of
$\hat{\alg{A}}$. Being a bijective homeomorphism and a homomorphism, the map
$\mu_{\alg{A}}$ is then an isomorphism of $\cate{M}_\tau$.
\end{proof}

The topology of the corollary is Roller's \emph{convex topology}, generated by
the prime convex subsets of $\alg{A}$, and the equivalence is his
\cite[Corollary~5.10]{Roller2016}. Since prime convex subsets that meet
pairwise have the finite intersection property, by Helly's theorem for median
algebras \cite[Theorem~2.2]{Roller2016}, this topology makes $\alg{A}$ compact
if and only if every family of pairwise intersecting prime convex subsets of
$\alg{A}$ has a nonempty intersection \cite[remark after
Proposition~5.9]{Roller2016}, that is, if and only if every morphism
$\alg{A}^+ \to \utilde{\alg{2}}$ is the evaluation at an element of $\alg{A}$.
Roller shows that this happens exactly when $\alg{A}^+$ contains no infinite
chain and no infinite family of pairwise transverse elements
\cite[Theorem~5.12]{Roller2016}.

\subsection{Profinite completion as natural extension}\label{subsec:natext}

The previous subsection built $\hat{\alg{A}}$ as the closure of $\alg{A}$ in
the product of all its finite quotients. This subsection shows that the
homomorphisms into $\alg{2}$ suffice, by identifying $\hat{\alg{A}}$ with the
closure of $\alg{A}$ in a power of the two element algebra, known as the
natural extension of $\alg{A}$. The identification is
once more delivered by the machinery of natural dualities recalled in
Subsection~\ref{subsec:nutshell}, and it comes in two steps. The alter ego that
produced the dualities of Theorem~\ref{thm:dualities} also computes the natural
extension, namely as the double dual $(\alg{A}^+)_+$
(Proposition~\ref{prop:natext}), while a general theorem on natural extensions
identifies that object with the profinite completion (Theorem~\ref{thm:DGHP}).
Read in succession, the two steps describe $\hat{\alg{A}}$ as the algebra of
ultrafilters on the bounded poc set $\alg{A}^+$
(Corollary~\ref{prop:profinite}).

That last description is the one the geometric literature arrives at by other
routes, and Subsection~\ref{subsec:roller} reads it off as the Roller
compactification. A single machine therefore yields the two dualities of
Section~\ref{sec:dualities}, the profinite completion, and a compactification
devised in geometric group theory.

The natural extension of $\alg{A}$ in $\cate{M}$, denoted
$n_\cate{M}(\alg{A})$, is constructed as follows. The evaluation map
\[
  e_{\alg{A}} \colon \alg{A} \longrightarrow \alg{2}_\tau^{\alg{A}^*},
  \qquad a \mapsto (h \mapsto h(a)),
\]
embeds $\alg{A}$ into the product $\alg{2}_\tau^{\alg{A}^*}$.
The \emph{natural extension} is the topological closure
\[
  n_\cate{M}(\alg{A}) := \overline{e_{\alg{A}}(A)} \;\subseteq\;
  \alg{2}_\tau^{\alg{A}^*},
\]
a closed subalgebra and hence an object of $\cate{M}_\tau$. Given a homomorphism $f \colon \alg{A} \to
\alg{B}$ in $\cate{M}$, the map $n_\cate{M}(f) \colon n_\cate{M}(\alg{A}) \to
n_\cate{M}(\alg{B})$ defined by
\[
  n_\cate{M}(f)(b)(h) := b(h \circ f), \qquad b \in n_\cate{M}(\alg{A}),
  \quad h \in \alg{B}^*,
\]
is a continuous homomorphism, making $n_\cate{M} \colon \cate{M} \to \cate{M}_\tau$
a functor \cite[Proposition~3.2]{DaveyGouveiaHaviarPriestley2011}.

\begin{proposition}\label{prop:natext}
 For any median algebra $\alg{A}$, the natural extension $n_\M(\alg{A})$ is
equal to $(\alg{A}^+)_+$.
\end{proposition}
\begin{proof}
The set $(\alg{A}^+)_+$ is a subalgebra of $\alg{2}_\tau^{\alg{A}^+}$
(Definition~\ref{def:complex}), closed because the conditions defining a
morphism
$\alg{A}^+ \to \utilde{\alg{2}}$ involve at most two coordinates at a time. It
contains $e_{\alg{A}}(A)$, so it suffices to prove the inclusion
$(\alg{A}^+)_+ \subseteq \overline{e_{\alg{A}}(A)}$.
Let $f\in (\alg{A}^+)_+$ and $F$ be a finite subset of $\alg{A}^*$. The
restriction of $f$ to the finite substructure of $\alg{A}^*$ formed by $F$,
$F^\co$ and the two constant maps is a morphism into
$\utilde{\alg{2}}_\tau$. By injectivity of $\utilde{\alg{2}}_\tau$ in
$\Poc_\tau$, it extends to an element of $(\alg{A}^*)_*$, which is
$e_{\alg{A}}(a)$ for some $a \in A$, hence $a$ witnesses $f$ on $F$.
\end{proof}

This first step is the median instance of a general fact. If an alter ego
$\utilde{\alg{N}}$ yields a duality on $\mathcal{A}$ at the finite level, then
the natural extension of every $\alg{A} \in \mathcal{A}$ consists of all the
maps $D(\alg{A}) \to N$ that preserve the structure of $\utilde{\alg{N}}$,
continuous or not \cite[Theorem~2.3]{DaveyHaviarPriestley2012}. The duality at
hand is strong by Theorem~\ref{thm:dualities}(\ref{it:dual-poctau}), hence a
duality at the finite level, and the description it yields is the composite
$-_+ \circ -^\flat \circ -^*$ of the opening of this section, whose value at
$\alg{A}$ is the double dual $(\alg{A}^+)_+$.

The second step is the specialization to median algebras of
\cite[Theorem~3.6]{DaveyGouveiaHaviarPriestley2011}, where the isomorphism is
obtained abstractly, from the functoriality of the natural extension and the
universal property of the limit. We give an explicit proof because the map
admits a direct description. The rest of the argument is the same in both
proofs.

\begin{theorem}\label{thm:DGHP}
For every median algebra $\alg{A}$, there is an isomorphism 
\[
	n_\M(\alg{A}) \;\cong\; \hat{\alg{A}}
\]
of Boolean topological algebras that fixes $\alg{A}$ and is natural in
$\alg{A}$.
\end{theorem}

\begin{proof}
We give an explicit construction of an
isomorphism $v_{\alg{A}}\colon n_{\M}(\alg{A}) \to \hat{\alg{A}}$. Let $f\in
n_{\M}(\alg{A})$ and $\alpha \in \mathcal{S}_{\alg{A}}$. The structure
$(\alg{A}/\alpha)^*$ embeds as a finite discrete substructure $F_\alpha$ of
$\alg{A}^*$ via $(\pi_{\alpha})^*$, where $\pi_\alpha$ is the quotient map
$\alg{A} \to\alg{A}/\alpha$. Since $f$ is locally in $\alg{A}$, there is a
witness $a\in A$ for $f$ on $F_\alpha$, so that
$f|_{F_\alpha}=e_{\alg{A}}(a)|_{F_\alpha}$. We set
\[
	v_{\alg{A}}(f)_\alpha := a/\alpha.
\]
\textit{Well-definedness.} Assume that $a$ and $a'$ both witness $f$ on $F_\alpha$.
Since $F_\alpha=\{u\circ \pi_{\alpha} \mid u \in (\alg{A}/\alpha)^*\}$, we obtain
$e_{\alg{A}}(a)(u \circ \pi_\alpha) = e_{\alg{A}}(a')(u \circ \pi_\alpha)$ for
every $u \in (\alg{A}/\alpha)^*$. It follows that $a/\alpha = a'/\alpha$ since
$(\alg{A}/\alpha)^*$ separates the elements of $\alg{A}/\alpha$.  

\textit{Thread condition.} Let $G$ be a finite subset of
$\mathcal{S}_{\alg{A}}$, with meet $\alpha$, and let $a$ witness $f$ on
$F_\alpha$. For every $\beta \in G$ we have $F_\beta \subseteq F_\alpha$, since
every $u \in (\alg{A}/\beta)^*$ satisfies $u \circ \pi_\beta = (u \circ
\varphi_{\alpha\beta}) \circ \pi_\alpha$. Hence $a$ also witnesses $f$ on
$F_\beta$, so that $v_{\alg{A}}(f)_\beta = a/\beta$ for every $\beta \in G$.
This shows that $v_{\alg{A}}(f)$ belongs to $\hat{\alg{A}}$.

\textit{Homomorphism.} Let $f, g, h \in
n_{\M}(\alg{A})$ and $\alpha \in \mathcal{S}_{\alg{A}}$, and choose witnesses
$a, b, c \in A$ for $f, g, h$ on $F_\alpha$. Since the median operation
of $\alg{2}_\tau^{\alg{A}^*}$ is computed coordinatewise, the element
$m(a,b,c)$ witnesses $m(f,g,h)$ on $F_\alpha$, whence
\[
  v_{\alg{A}}(m(f,g,h))_\alpha = m(a,b,c)/\alpha
  = m\bigl(v_{\alg{A}}(f), v_{\alg{A}}(g), v_{\alg{A}}(h)\bigr)_\alpha.
\]
Moreover, for every $a \in A$ the element $a$ witnesses $e_{\alg{A}}(a)$ on
every $F_\alpha$, so that
\begin{equation}\label{eq:v-fixes-A}
  v_{\alg{A}} \circ e_{\alg{A}} = \mu_{\alg{A}}
\end{equation}
which shows that the map $v_{\alg{A}}$ fixes $\alg{A}$.

\textit{Continuity.} The topology of $\hat{\alg{A}}$ is the initial topology
with respect to the projections $p_\alpha \colon \hat{\alg{A}} \to
\alg{A}/\alpha$, so it suffices to prove that each $p_\alpha \circ
v_{\alg{A}}$ is continuous. If $f$ and $g$ agree on $F_\alpha$, then any
witness for $f$ on $F_\alpha$ also witnesses $g$, so that $v_{\alg{A}}(f)_\alpha
= v_{\alg{A}}(g)_\alpha$. Hence $p_\alpha \circ v_{\alg{A}}$ factors through
the restriction map $n_{\M}(\alg{A}) \to \alg{2}^{F_\alpha}$, which is
continuous by definition of the product topology and has a finite discrete
codomain.

\textit{Injectivity.} Let $f$ and $g$ be distinct elements of
$n_{\M}(\alg{A})$, and let $h \in \alg{A}^*$ satisfy $f(h) \neq g(h)$. The
kernel $\alpha$ of $h$ belongs to $\mathcal{S}_{\alg{A}}$, and $h$ factors as
$u \circ \pi_\alpha$ for some $u \in (\alg{A}/\alpha)^*$, so $h \in F_\alpha$.
If $a$ and $b$ witness $f$ and $g$ on $F_\alpha$, respectively, then
\[
  u(a/\alpha) = h(a) = f(h) \neq g(h) = h(b) = u(b/\alpha),
\]
whence $v_{\alg{A}}(f)_\alpha \neq v_{\alg{A}}(g)_\alpha$.

\textit{Surjectivity.} 
The algebra $n_{\M}(\alg{A})$ is compact and the map
$v_{\alg{A}}$ is continuous, so the image of $v_{\alg{A}}$ is a compact, hence
closed, subset of the Hausdorff space $\hat{\alg{A}}$. By
\eqref{eq:v-fixes-A}, this image contains $\mu_{\alg{A}}(A)$, which is dense in
$\hat{\alg{A}}$. A closed dense subset is the whole space, so $v_{\alg{A}}$ is
surjective. A continuous bijection from a compact space to a Hausdorff space is
a homeomorphism, whence $v_{\alg{A}}$ is an isomorphism in $\cate{M}_\tau$.

\textit{Naturality.} Let $f \colon \alg{A} \to \alg{B}$ be a homomorphism in
$\cate{M}$. We have to show that $\hat{f} \circ v_{\alg{A}} = v_{\alg{B}} \circ
n_{\M}(f)$. Both sides are continuous maps from $n_{\M}(\alg{A})$ into the
Hausdorff space $\hat{\alg{B}}$, so by the opening remark of this section it
suffices to verify equality on the dense subset $e_{\alg{A}}(A)$. For $a \in A$
and $h \in \alg{B}^*$ we have $n_{\M}(f)(e_{\alg{A}}(a))(h) = e_{\alg{A}}(a)(h
\circ f) = h(f(a))$, that is, the identity $n_{\M}(f) \circ e_{\alg{A}} =
e_{\alg{B}} \circ f$ holds. Combining \eqref{eq:v-fixes-A} with the naturality
identity \eqref{eq:mu-natural}, we obtain
\[
  v_{\alg{B}}\bigl(n_{\M}(f)(e_{\alg{A}}(a))\bigr)
  = v_{\alg{B}}(e_{\alg{B}}(f(a)))
  = \mu_{\alg{B}}(f(a))
  = \hat{f}(\mu_{\alg{A}}(a))
  = \hat{f}\bigl(v_{\alg{A}}(e_{\alg{A}}(a))\bigr),
\]
which concludes the proof.
\end{proof}

\begin{corollary}\label{prop:profinite}
For every median algebra $\alg{A}$, there is an isomorphism
\[
  (\alg{A}^+)_+ \;\cong\; \hat{\alg{A}}
\]
of Boolean topological median algebras that fixes $\alg{A}$ and is natural in
$\alg{A}$.
\end{corollary}

We will use the following preservation lemma in the next section.

\begin{lemma}\label{lem:pres}
Let $s(\mathbf{x}, \mathbf{y})$ and $t(\mathbf{x}, \mathbf{y})$ be terms of
median algebras, and let $\alg{A}$ be a median algebra. For any $B \subseteq A$
and $\mathbf{c}, \mathbf{d} \in A^{|\mathbf{y}|}$, if $t(\mathbf{a},
\mathbf{c}) = s(\mathbf{a}, \mathbf{d})$ for every $\mathbf{a} \in
B^{|\mathbf{x}|}$, then $t(\mathbf{a}, \mathbf{c}) = s(\mathbf{a},
\mathbf{d})$ for every $\mathbf{a} \in \overline{B}^{|\mathbf{x}|}$, where
$\overline{B}$ denotes the closure of $B$ in $\hat{\alg{A}}$.
\end{lemma}
\begin{proof}
Let $\mathbf{a} \in \overline{B}^{|\mathbf{x}|}$ and $u \in \alg{A}^*$. For each
element $a_i$, locality provides a witness $a'_i \in B$ satisfying $a'_i(u) = a_i(u)$. Since
median algebra terms are evaluated coordinatewise in $\hat{\alg{A}} \subseteq
\alg{2}^{\alg{A}^*}$, we obtain
\[
t(\mathbf{a}, \mathbf{c})(u) = t(\mathbf{a}', \mathbf{c})(u)
= s(\mathbf{a}', \mathbf{d})(u) = s(\mathbf{a}, \mathbf{d})(u),
\]
where the middle equality uses $\mathbf{a}' \in B^{|\mathbf{x}|}$ and the
hypothesis. Since $u$ was arbitrary, the identity $t(\mathbf{a}, \mathbf{c}) =
s(\mathbf{a}, \mathbf{d})$ holds.
\end{proof}

\subsection{The Roller compactification}\label{subsec:roller}

Read through the dictionary of Remark~\ref{rem:terminology},
Corollary~\ref{prop:profinite} recovers a construction of geometric group
theory.

Let $\alg{A}$ be the median algebra of the $0$-skeleton of a $\mathrm{CAT}(0)$
cube complex $X$. The bounded poc set $\alg{A}^+$ is the halfspace system
$\mathcal{H}(X)$, and $(\alg{A}^+)_+$ is the median algebra of ultrafilters on
$\mathcal{H}(X)$, topologized as a closed subspace of
$\{0,1\}^{\mathcal{H}(X)}$. This is the \emph{Roller compactification}
$\overline{X}$ of $X$ \cite{Roller2016}, so that
Corollary~\ref{prop:profinite} reads
\[
  \overline{X} \;\cong\; \hat{\alg{A}}.
\]
For an arbitrary median algebra $\alg{A}$, the double dual $(\alg{A}^+)_+$ is
Bowditch's algebra $\mathcal{F}(\mathcal{H}(\alg{A}))$ of flows on the
halfspaces of $\alg{A}$, the complement of $\alg{A}$ in it being the
\emph{Roller boundary} $\partial_R \alg{A}$ \cite[\S11.12]{Bowditch2022}, and
it is Fioravanti's double dual $\alg{A}^{\circ\circ}$
\cite[Definition~4.2 and Lemma~4.3(1)]{Fioravanti2020}. On objects, the pair
$-^+$ and $-_+$ is the correspondence between halfspaces and flows of
Bowditch \cite[\S9.2]{Bowditch2022}.

A second construction of that literature must be distinguished from it,
because it produces a different and in general smaller object. The
\emph{zero-completion} of $\alg{A}$ was introduced by Bandelt and Meletiou as
the algebra of compatible semilattice operations on $\alg{A}$
\cite[Theorem~1]{BandeltMeletiou1993}, and identified by them with the largest
subalgebra of $\alg{2}^X$ in which $\alg{A}$ is convex, whenever $\alg{A}$ is a
subalgebra of $\alg{2}^X$ with convex hull $\alg{2}^X$
\cite[Theorem~2]{BandeltMeletiou1993}. Fioravanti realized it as the inverse
limit of the intervals of $\alg{A}$ under the so-called gate projections
$x \mapsto m(u,v,x)$ onto the intervals $[u;v]$
\cite[\S4.1]{Fioravanti2020}. It always embeds in $\hat{\alg{A}}$
\cite[Lemma~4.3(2)]{Fioravanti2020}, with equality when the intervals of
$\alg{A}$ are finite \cite[\S12.6]{Bowditch2022}, which covers the $0$-skeleta
above, but not in general \cite[Example~4.5]{Fioravanti2020}.
Remark~\ref{rem:ccan-vs-zero} below compares the two completions once the
interval-by-interval description of $\hat{\alg{A}}$ is available.

The definition of $\hat{\alg{A}}$ as a limit of finite quotients has a
geometric reading as well, since its finite quotients are the collapses of the
geometric literature. For a finite set $F$ of prime convex subsets of
$\alg{A}$, let $\theta_F$ be the congruence that identifies two elements of
$\alg{A}$ when no member of $F$ contains one of them without the other. It is
the kernel of the homomorphism $a \mapsto (h(a))_{h \in F}$ into
$\alg{2}^F$, hence of finite index. Every congruence $\alpha$ of finite index
arises in this way. Indeed, the finite quotient $\alg{A}/\alpha$ lies in $\M =
\classop{ISP}(\alg{2})$, so the finitely many homomorphisms $u \colon
\alg{A}/\alpha \to \alg{2}$ separate its elements. Their compositions $u \circ
\pi_\alpha$ with the quotient map form a finite set $F$, and $\theta_F$ is the
kernel of $\pi_\alpha$, which is $\alpha$. When $\alg{A}$ is the vertex set of a
median graph, the quotient
$\alg{A}/\theta_F$ is the \emph{hyperplane-collapse} of the graph along the
hyperplanes that bound the members of $F$, obtained by contracting every edge
dual to any other hyperplane \cite[Corollary~2.3.6]{Genevois2024}, and Caprace
and Sageev call it a restriction quotient \cite[\S2.3]{CapraceSageev2011}. Read
through this dictionary, Corollary~\ref{prop:profinite} says that the Roller
compactification is the limit of the finite hyperplane-collapses, a statement we
have not found in print, although its ingredients are standard.

This description is the instance, for median algebras, of the
coincidence of the natural extension with the profinite completion
(Theorem~\ref{thm:DGHP}), combined with the identification of the natural
extension with $(\alg{A}^+)_+$ (Proposition~\ref{prop:natext}). The same
apparatus delivers with it the compactness of the Roller compactification.

The universal property of Proposition~\ref{prop:ump-reflection}, combined
with the duality of Theorem~\ref{thm:dualities}(\ref{it:dual-poc}), describes
all the continuous images of $\hat{\alg{A}}$. A \emph{bounded sub-poc set} of
$\alg{A}^+$ is a subset $P$ that contains $\emptyset$ and $A$ and is closed
under complement, ordered by inclusion. It is \emph{separating} if it
separates the points of $\alg{A}$. The homomorphism $e_P \colon \alg{A} \to
P_+$ given by $a \mapsto \{h \in P \mid a \in h\}$ identifies two elements
that no member of $P$ separates, so it is injective exactly when $P$ is
separating.

\begin{proposition}\label{prop:compactifications}
Let $\alg{A}$ be a median algebra.
\begin{enumerate}
\item For every bounded sub-poc set $P$ of $\alg{A}^+$, the homomorphism
$e_P$ has dense image in $P_+$.
\item Let $f \colon \alg{A} \to \alg{B}$ be a homomorphism with dense image
into some $\alg{B} \in \M_\tau$. Then $P_f := \{v \circ f \mid v \in
F(\alg{B})\}$ is the only bounded sub-poc set $P$ of $\alg{A}^+$ for which an
isomorphism $\alg{B} \cong P_+$ carries $f$ to $e_P$. In particular, the
algebra $\alg{B}$ is a continuous image of $\hat{\alg{A}}$.
\item If the intervals of $\alg{A}$ are finite, then $\alg{A}^+$ is its only
separating bounded sub-poc set.
\end{enumerate}
\end{proposition}
\begin{proof}
(1) A basic neighbourhood of $x \in P_+$ is determined by finitely many
coordinates, which we may take to be members $u_1, \ldots, u_n$ of $P$ with
$x(u_i) = 1$, replacing a member by its complement if needed. As $x$ is a
morphism, no $u_i$ is contained in the complement of another one, so the
sets $u_1, \ldots, u_n$ meet pairwise. Finitely many convex subsets of a
median algebra that meet pairwise have a common element, by Helly's theorem
\cite[Theorem~2.2]{Roller2016}. So the sets $u_i$ have a common element $a$,
and $e_P(a)$ lies in the neighbourhood.

(2) Precomposition with $f$ maps $F(\alg{B})$ onto $P_f$, preserving the
structure, and it reflects the order. Indeed, if $v \circ f \leq w \circ f$,
then the clopen set $\{y \in B \mid v(y) > w(y)\}$ misses the dense set
$f(A)$, hence is empty, and $v \leq w$. So $P_f$ is a bounded sub-poc set of
$\alg{A}^+$ isomorphic to $F(\alg{B})$, and
Theorem~\ref{thm:dualities}(\ref{it:dual-poc}) yields an isomorphism
$\alg{B} \cong F(\alg{B})_+ \cong (P_f)_+$, which sends $f(a)$ to the
morphism $v \circ f \mapsto v(f(a))$, that is, to $e_{P_f}(a)$. Conversely,
let $\varphi \colon \alg{B} \to P_+$ be an isomorphism with $\varphi \circ f =
e_P$. By the same theorem, the continuous homomorphisms $P_+ \to \alg{2}$ are
the evaluations at the members of $P$, so $F(\alg{B})$ consists of the maps
$y \mapsto \varphi(y)(u)$ for $u \in P$. Precomposing them with $f$ yields
the maps $a \mapsto e_P(a)(u) = u(a)$, whence $P_f = P$. Finally, the
extension $\hat{\alg{A}} \to \alg{B}$ of $f$ given by
Proposition~\ref{prop:ump-reflection} is onto, its image being compact and
dense.

(3) A proper prime convex subset $h$ of $\alg{A}$ separates some $x$ and $y$
with $[x;y] = \{x, y\}$, and $h$ is then the only prime convex subset
containing $x$ but not $y$, as a point $z$ in exactly one of two such subsets
would make $m(x, y, z)$ equal to both $x$ and $y$.
\end{proof}

\begin{remark}\label{rem:compactifications}
Proposition~\ref{prop:compactifications} is the median counterpart of the
classical description of the zero-dimensional compactifications of a space
by separating families of clopen sets, and differs from
\cite[Remark~5.7]{Roller2016}, by which a median algebra carries at most one
Boolean topology. Separating sub-poc sets give the compactifications of
$\alg{A}$, and the others those of its collapses. Geometrically, the Roller
compactification of a $\mathrm{CAT}(0)$ cube complex is its only
zero-dimensional median compactification. That the vertex set is dense and
convex in it is in \cite[Theorem~C]{Fioravanti2020}, but not this
uniqueness. Without finite intervals there are many. In the real line, the
rational cuts form a countable separating sub-poc set, whose complex algebra
is metrizable and adds two points beside each rational number and two
points at infinity, whereas $\hat{\R}$ adds two points beside every real
number and two points at infinity. A group acting on
$\alg{A}$ permutes the bounded sub-poc sets of $\alg{A}^+$, and each
invariant one yields an equivariant compactification.
\end{remark}

\begin{remark}\label{rem:recognition}
Let $Y$ be a component of the Roller compactification of a median graph. It
has finite intervals, and its closure in the Roller compactification is a
Boolean topological median algebra in which $Y$ is dense
(Proposition~\ref{prop:profinite-char}). By
Proposition~\ref{prop:compactifications}, this closure is the Roller
compactification of $Y$, embedded by the unique continuous extension of the
inclusion (Proposition~\ref{prop:ump-reflection}). This recovers
\cite[Proposition~2.2.1 and Exercise~2.5.10(3)]{Genevois2024}.
\end{remark}

\section{Profinite completion through canonical
extensions}\label{sec:can-ext}

J\'onsson and Tarski introduced canonical extensions to extend Stone's
representation theorem to Boolean algebras with operators
\cite{JonssonTarski1951,JonssonTarski1952}.
They have since become an algebraic tool for completeness and correspondence
results in modal and substructural logic \cite{Gehrke2004,Dunn2005,Gehrke2014},
and Fussner and Gehrke give a systematic account \cite{FussnerGehrke2025}. The
relation of canonical extensions to the
profinite completion is a recurring theme: in a finitely generated variety of
lattice-based algebras, the two constructions coincide
\cite{Harding2006,Gouveia2010,Vosmaer2010}.  Canonical extensions have been
further constructed for general posets \cite{Dunn2005}, and the case of
semilattices has received a dedicated treatment \cite{Gouveia2013}.  Having
built $\hat{\alg{A}}$ in Section~\ref{sec:profinite}, we ask whether it is a
canonical extension of any of the median semilattices $(\alg{A}, \leq_a)$ of
$\alg{A}$.

When $(\alg{A}, \leq_a)$ has a top, the whole of $(\hat{\alg{A}}, \leq_a)$
is its canonical extension (Corollary~\ref{cor:global-top}). A median
semilattice has in general no top, however, and then ${(\hat{\alg{A}}, \leq_a)}$
need not be a complete lattice, as for the chain $\Z$ with basepoint
$0$. Interval by interval, the profinite completion still realizes canonical
extensions (Theorem~\ref{thm:can-ext}), and this interval-wise property
characterizes the profinite completion among the dense extensions of $\alg{A}$
in $\M_\tau$ (Theorem~\ref{thm:char-profinite}).

Order completions can be approached from three points of view
\cite{FussnerGehrke2025}. There is a concrete construction (obtained through
duality), an internal characterization (saying how the original structure sits
inside the completion), and a universal property (presenting the completion as
a left adjoint to a forgetful functor). For the profinite completion the first
and the third are those of Section~\ref{sec:profinite}
(Corollary~\ref{prop:profinite} and Proposition~\ref{prop:ump-reflection}),
and this section supplies the second (Theorem~\ref{thm:char-profinite}).

Median semilattices need not be lattices, but each of their finite upper-bounded
subsets has a least upper bound. The following definition fixes the vocabulary
for a weaker form of order completeness, and for the sub-semilattices that
keep
track of the existing joins.

\begin{definition}[\cite{Davey2002}]\label{def:cond-complete}
A sub-semilattice $\alg{S}'$ of a
meet-semilattice $\alg{S}$ is a \emph{subalgebra} of $\alg{S}$ if $\alg{S}'$ is
closed under all existing nonempty finite joins in $\alg{S}$. It is a \emph{complete subalgebra} of
$\alg{S}$ if $\alg{S}'$ is closed under all existing nonempty joins in $\alg{S}$.

A meet-semilattice $\alg{S}$ is \emph{conditionally complete} if every
upper-bounded subset of $S$, the empty subset included, has a least upper
bound.
\end{definition}

\begin{remark}\label{rem:comp}
\begin{enumerate}
\item\label{it:mdf}
Complete lattices are instances of conditionally complete meet-semilattices.
Every conditionally complete meet-semilattice has a bottom element
$\bigvee \varnothing$. A meet-semilattice $(S, \wedge)$ is conditionally
complete if and only if every \emph{nonempty} subset $S'$ of $S$ has a
greatest
lower bound $\bigwedge S'$.
\item A conditionally complete meet-semilattice with a top element is a
complete lattice.
\item The bottom element is not part of the signature. Accordingly, a
subalgebra of a meet-semilattice that has a bottom element need not contain
it, and a homomorphism of meet-semilattices need not preserve
bottom elements. The next item provides a complete subalgebra that omits
the bottom.
\item If $\alg{S}'$ is a complete subalgebra of a meet-semilattice $\alg{S}$, then $\alg{S}'$ need
not be closed under all nonempty meets of  $\alg{S}$. For instance, the cofinite
subsets of $\N$ form a complete subalgebra of $(2^\N, \wedge)$, while
$\bigwedge \{\N\setminus\{n\} \mid n \in \N\}$ exists in $2^\N$ and is equal to
$\varnothing$, which is not cofinite.
\end{enumerate}
\end{remark}

Following standard terminology, a \emph{filter} of a poset $\alg{P}$ is a
nonempty down-directed up-set of $\alg{P}$, and an \emph{ideal} of $\alg{P}$ is
a nonempty up-directed down-set of $\alg{P}$. 

\begin{definition}[\cite{Dunn2005}]\label{def:can-ext}
If $\alg{P}$ is a poset and $\alg{L}$ is a complete lattice such that
$\alg{P}$ is a subposet of $\alg{L}$, then $\alg{L}$ is called a
\emph{completion} of $\alg{P}$. Under this condition, an element $x$ of
$\alg{L}$ is a \emph{filter element} (or \emph{closed} element) if $x =
\bigwedge F$ for some filter $F$ of $\alg{P}$, and an \emph{ideal element}
(or \emph{open} element) if $x = \bigvee I$ for some ideal $I$ of $\alg{P}$.

A completion $\alg{L}$ of $\alg{P}$ is
\begin{enumerate}
\item \emph{dense} if every element of $\alg{L}$ is both a join of filter elements and a meet of ideal elements,
\item \emph{compact} if for every filter $F$ and every ideal $I$ of $\alg{P}$
such that $\bigwedge F \leq \bigvee I$, there are $s\in F$ and $t\in I$ such
that $s \leq t$, or equivalently $F \cap I \neq \varnothing$.
\end{enumerate}
A completion of $\alg{P}$ is called a \emph{canonical extension} if it is both dense and compact.
\end{definition}

Every poset $\alg{P}$ has a unique canonical extension, up to an isomorphism
fixing $\alg{P}$ \cite{Dunn2005}.

The complete lattice structures obtained in this section all stem from the
following topological observation.

\begin{lemma}\label{lem:closed-joins}
Let $X$ be a set and let $x \in \alg{2}^X$.
\begin{enumerate}
\item\label{it:dir-loc} The join in $(\alg{2}^X, \leq_x)$ of any up-directed
subset $D$ of $\alg{2}^X$ belongs to the closure of $D$ in
$\alg{2}_\tau^X$.
\item\label{it:closed-joins} If a closed subset $C$ of $\alg{2}_\tau^X$ is
closed under nonempty finite joins in $(\alg{2}^X, \leq_x)$, then $C$ is
closed under arbitrary nonempty joins.
\end{enumerate}
The dual statements hold for meets.
\end{lemma}
\begin{proof}
\eqref{it:dir-loc} Let $s$ be the join of $D$, computed coordinatewise, and let
$F \Subset X$. For each $u \in F$, the value $s(u)$ is the join
of the subset $\{t(u) \mid t \in D\}$ of $(\{0,1\}, \leq_{x(u)})$, so there is
$t_u \in D$ with $t_u(u) = s(u)$.  Since $D$ is up-directed, some $t^* \in D$ is
above every $t_u$. Then $t^*$ agrees with $s$ on $F$, since $s(u) \geq_{x(u)}
t^*(u) \geq_{x(u)} t_u(u) = s(u)$ for every $u \in F$. Hence $s$ is locally in
$D$.

\eqref{it:closed-joins} Let $S$ be a nonempty subset of $C$, and let $D$ be
the set of joins of the nonempty finite subsets of $S$. Then $D$ is an
up-directed subset of $C$, and the joins of $S$ and of $D$ coincide. By
\eqref{it:dir-loc}, this join belongs to $\overline{D} \subseteq C$.
\end{proof}

\begin{lemma}\label{lem:ont}
For any $\alg{A} \in \M_\tau$ and any $a,b \in A$, the interval $[a;b]$ is a
complete and completely distributive bounded lattice with respect to $\leq_a$
and a complete subalgebra of $(\alg{A}, \leq_a)$.
Its nonempty joins and meets are computed coordinatewise in any power of
$\alg{2}_\tau$ in which $\alg{A}$ embeds.
\end{lemma}

\begin{proof}
Assume that $\alg{A}$ embeds in $\alg{2}_\tau^X$ both algebraically and
topologically. The interval $[a;b]_{\alg{A}}$ is a closed subspace of
$\alg{2}_\tau^X$, being the equalizer of two term functions. It is closed under
finite meets and joins, since they are computed term-wise.
Lemma~\ref{lem:closed-joins}\eqref{it:closed-joins} shows that
$[a;b]_{\alg{A}}$ is closed under all nonempty joins and meets of
$(\alg{2}^X, \leq_a)$, hence a completely distributive complete lattice with
bounds $a$ and $b$.

To see that $[a;b]_{\alg{A}}$ is also a complete subalgebra of $(\alg{A},
\leq_a)$, let $S$ be a nonempty subset of it having a join in $(\alg{A},
\leq_a)$. That join coincides with the join of $S$ computed in $(\alg{2}^X,
\leq_a)$, since the latter lies in $[a;b]_{\alg{A}}$ by the preceding paragraph
and is below every upper bound of $S$ in $(\alg{A}, \leq_a)$.
\end{proof}

The following lemma shows how the profinite completion produces complete
lattices from bounded subalgebras.

\begin{lemma}\label{lem:closure}
Let $\alg{A} \in \M$, let $a \in A$, and let $P$ be a bounded subalgebra
of $(\alg{A}, \leq_a)$ with bottom $c$ and top $d$. The closure $\overline{P}$
of $P$ in $\hat{\alg{A}}$ is a complete lattice under $\leq_a$, with bottom
$c$ and top $d$, and is closed under all nonempty joins and meets of
$(\alg{2}^{\alg{A}^*}, \leq_a)$.
\end{lemma}
\begin{proof}
By Lemma~\ref{lem:pres} applied to the equations $m(a, c, x) = c$ and
$m(a, x, d) = x$, the elements $c$ and $d$ remain the bottom and top of
$\overline{P}$ in $(\hat{\alg{A}}, \leq_a)$. The binary meet of
$(\alg{2}^{\alg{A}^*}, \leq_a)$ is the term operation $m(a, -, -)$, and on the
elements below $d$ its binary join is $m(-, -, d)$ by
Lemma~\ref{lem:interval-dual}.
Being a sub-semilattice of $(\alg{A}, \leq_a)$ closed under its existing
joins, $P$ is closed under both. The same holds for $\overline{P}$, which lies
below $d$. Indeed, given $x, y \in \overline{P}$ and a finite $F \subseteq
\alg{A}^*$, witnesses $x', y' \in P$ for $x$ and $y$ on $F$ give elements $m(a,
x', y')$ and $m(x', y', d)$ of $P$ agreeing on $F$ with $m(a, x, y)$ and $m(x,
y, d)$, since terms are evaluated coordinatewise.

The closure $\overline{P}$ is a closed subset of
$\alg{2}_\tau^{\alg{A}^*}$, so Lemma~\ref{lem:closed-joins}\eqref{it:closed-joins} shows that
$\overline{P}$ is closed under all nonempty joins and meets of
$(\alg{2}^{\alg{A}^*}, \leq_a)$. Since $\overline{P}$ is also bounded, it is a
complete lattice.
\end{proof}

\begin{theorem}\label{thm:can-ext}
Let $\alg{A} \in \M$, let $a \in A$, and let $P$ be a bounded subalgebra of
$(\alg{A}, \leq_a)$. The closure $\overline{P}$ of $P$ in
$\hat{\alg{A}}$ is a canonical extension of $P$.
\end{theorem}
\begin{proof}
By Lemma~\ref{lem:closure}, the closure $\overline{P}$ is a completion of
$P$. It remains to prove that this completion is compact and dense.

\textit{Compactness.}
Suppose $\mathcal{F}$ is a filter and $\mathcal{I}$ an ideal of $(P, \leq_a)$
with $s \not\leq_a t$ for every $s \in \mathcal{F}$ and $t \in \mathcal{I}$. For
each such pair, set
\[
  C_{s,t} := \{h \in \alg{A}^* \mid s(h) >_{a(h)} t(h)\}.
\]
The family $\{C_{s,t}\mid s \in \mathcal{F}, t\in\mathcal{I}\}$ is composed of
nonempty closed subspaces of $\alg{A}^*$ and has the finite intersection
property, since for $s_0 \leq_a s_1, s_2$ in $\mathcal{F}$ and
$t_0 \geq_a t_1, t_2$ in $\mathcal{I}$, one has $C_{s_0,t_0} \subseteq
C_{s_1,t_1} \cap C_{s_2,t_2}$. By compactness of $\alg{A}^*$, one can pick $h \in
\bigcap_{s,t} C_{s,t}$, so
\[
s(h) >_{a(h)} t(h), \qquad \text{for all } s\in \mathcal{F} \text{ and } t \in \mathcal{I}.
\]
By Lemma~\ref{lem:closure}, $\bigvee_{\overline{P}} \mathcal{I} =
\bigvee_{\alg{2}^{\alg{A}^*}} \mathcal{I}$ and $\bigwedge_{\overline{P}}
\mathcal{F} =
\bigwedge_{\alg{2}^{\alg{A}^*}} \mathcal{F}$. Hence
\[
\Bigl(\bigvee\nolimits_{\overline{P}} \mathcal{I}\Bigr)(h)
<_{a(h)}
\Bigl(\bigwedge\nolimits_{\overline{P}} \mathcal{F}\Bigr)(h),
\]
and therefore $\bigwedge_{\overline{P}} \mathcal{F} \not\leq_a
\bigvee_{\overline{P}} \mathcal{I}$.

\textit{Density.}
Let $x \in \overline{P}$. For each finite $F \subseteq \alg{A}^*$, let $G_F = \{c \in
P \mid c|_F = x|_F\}$ be the set of witnesses for $x$ on $F$ in $P$; it is
nonempty since $x$ is locally in $P$. Let $K_F$ and
$I_F$ be the filter and ideal generated by $G_F$ in $(P, \leq_a)$. We claim that
\[
  x = \bigvee\{\bigwedge K_F \mid F \Subset \alg{A}^*\}
    = \bigwedge\{\bigvee I_F \mid F \Subset \alg{A}^*\},
\]
where the complete operations are computed relative to $\leq_a$.

For the first identity, let $h \in \alg{A}^*$. If $h \in F$, then every $c \in K_F$ satisfies
$c(h) \geq_{a(h)} x(h)$, since $c \geq_a g$ for some $g \in G_F$, and $g(h) =
x(h)$ as $h \in F$. Moreover, since $\emptyset\neq G_F \subseteq K_F$,
we obtain $(\bigwedge K_F)(h) = x(h)$.

If $h \notin F$, locality provides $c \in G_{F \cup \{h\}} \subseteq K_F$ with
$c(h) = x(h)$, so $(\bigwedge K_F)(h) \leq_{a(h)} c(h) = x(h)$.

Thus $\bigwedge K_F \leq_a x$ for every finite subset $F$ of
$\alg{A}^*$. Now, let $h \in \alg{A}^*$. Considering $F = \{h\}$ yields
$(\bigwedge K_{\{h\}})(h) = x(h)$, which proves the first identity.

The second identity follows by a dual argument.
\end{proof}

\begin{corollary}\label{cor:global-top}
Let $\alg{A} \in \M$ and $a \in A$. If the poset $(\alg{A}, \leq_a)$ has a
greatest element, then $(\hat{\alg{A}}, \leq_a)$ is the canonical extension
of $(\alg{A}, \leq_a)$. 
\end{corollary}

The interval case of Theorem~\ref{thm:can-ext} in fact characterizes the
profinite completion among the dense extensions of $\alg{A}$ in
$\cate{M}_\tau$.

\begin{theorem}\label{thm:char-profinite}
Let $\alg{A}\in \M$ be a median algebra and $\alg{A}'\in
\M_\tau$ an extension of $\alg{A}$ in which $\alg{A}$ is dense.
Then $\alg{A}'$ is a profinite completion of
$\alg{A}$ if and only if for every $a, b \in A$, the interval $[a;b]$ in
$(\alg{A}', \leq_a)$ is a canonical extension of $[a;b]$ in $(\alg{A}, \leq_a)$.
\end{theorem}
\begin{proof}
\textit{$(\Rightarrow)$}
Assume that $\alg{A}' = \hat{\alg{A}}$, and let $a, b \in A$. We first check
that the interval $[a;b]$ in $\hat{\alg{A}}$ is the closure of
$[a;b]_{\alg{A}}$. Intervals are closed, so the closure of $[a;b]_{\alg{A}}$
is contained in $[a;b]_{\hat{\alg{A}}}$. Conversely, if $z \in
[a;b]_{\hat{\alg{A}}}$, if $F$ is a finite subset of $\alg{A}^*$ and if $w
\in A$ witnesses $z$ on $F$, then the element $m(a, w, b)$ of
$[a;b]_{\alg{A}}$ also witnesses $z$ on $F$, since $z|_F = m(a|_F, z|_F,
b|_F) = m(a, w, b)|_F$. Hence the interval $[a;b]_{\hat{\alg{A}}}$ is the
closure of $[a;b]_{\alg{A}}$, which is a bounded subalgebra of $(\alg{A},
\leq_a)$, and Theorem~\ref{thm:can-ext} applied to $P = [a;b]_{\alg{A}}$
yields the conclusion.

\textit{$(\Leftarrow)$}
Assume $\alg{A}$ is dense in $\alg{A}'$ and every interval $[a;b]_{\alg{A}'}$
is a canonical extension of $[a;b]_{\alg{A}}$.
By Proposition~\ref{prop:ump-reflection}, the embedding $\alg{A}
\hookrightarrow \alg{A}'$ extends to a unique continuous median homomorphism
$\phi\colon \hat{\alg{A}} \to \alg{A}'$ fixing $\alg{A}$.

\textit{Injectivity on intervals.}
Let $u, v \in A$. If $z \in [u;v]_{\hat{\alg{A}}}$, then $\phi(z) =
\phi(m(u, z, v)) = m(u, \phi(z), v)$, so $\phi$ restricts to a continuous
median homomorphism $\psi \colon [u;v]_{\hat{\alg{A}}} \to
[u;v]_{\alg{A}'}$ fixing $[u;v]_{\alg{A}}$. Both intervals are canonical
extensions of $[u;v]_{\alg{A}}$, the former by the first part of the proof
and the latter by hypothesis. We show that $\psi$ is injective.

First, we prove that $\psi$ preserves filter and ideal elements. The term
functions $m(u, -, -)$ and $m(-, -, v)$ compute the
binary meet and join of $\leq_u$ on both intervals, the join by
Lemma~\ref{lem:interval-dual}. Since $\psi$ preserves them, it is
order-preserving. Let $F$ be a filter of $([u;v]_{\alg{A}}, \leq_u)$, and let
$k$ and $k'$ be the meets of $F$ computed in $[u;v]_{\hat{\alg{A}}}$ and in
$[u;v]_{\alg{A}'}$, respectively. The element $\psi(k)$ is a lower bound of $F$,
so $\psi(k) \leq k'$. Moreover, the meet $k$ is computed coordinatewise by
Lemma~\ref{lem:ont}, so the dual statement of
Lemma~\ref{lem:closed-joins}\eqref{it:dir-loc} places $k$ in the closure of $F$,
and by continuity $\psi(k)$ lies in the closure of $F$ computed in $\alg{A}'$.
Suppose $\psi(k) \neq k'$. Since $\alg{A}' \in \M_\tau$ is a closed subalgebra
of a power of $\alg{2}_\tau$, some continuous median homomorphism $\chi \colon
\alg{A}' \to \alg{2}_\tau$ separates $\psi(k)$ from $k'$. If $\chi(u) =
\chi(v)$, then $\chi(z) = \chi(m(u, z, v)) = \chi(u)$ for every $z \in
[u;v]_{\alg{A}'}$, and $\chi$ separates nothing in the interval. So, replacing
$\chi$ by its complement if necessary, we may assume $\chi(u) = 0$ and
$\chi(v) = 1$, so that $\chi$ preserves the binary meets and joins of $\leq_u$
on $[u;v]_{\alg{A}'}$ and is order-preserving there. From $\psi(k) \leq k'$ we
obtain $\chi(\psi(k)) = 0$ and $\chi(k') = 1$. Since $k' \leq s$ for every $s
\in F$, the map $\chi$ is constantly $1$ on $F$, hence on its closure, which
contains $\psi(k)$, a contradiction. Hence $\psi(k) = k'$, and $\psi$ preserves
filter elements. The dual argument applies to ideal elements.

Now let $x, y \in [u;v]_{\hat{\alg{A}}}$ with $x \neq y$, say $x \not\leq_u y$.
Since $x$ is a join of filter elements, some filter element $k = \bigwedge F$
satisfies $k \leq x$ and $k \not\leq y$. Since $y$ is a meet of ideal elements,
some ideal element $o = \bigvee I$ satisfies $y \leq o$ and $k \not\leq o$.
There are then no $s \in F$ and $t \in I$ with $s \leq t$, as this would give $k
\leq s \leq t \leq o$. By compactness of the canonical extension
$[u;v]_{\alg{A}'}$, the meet of $F$ is not below the join of $I$ in
$[u;v]_{\alg{A}'}$ either. Since $\psi(x) \geq \psi(k) = \bigwedge F$ and
$\psi(y) \leq \psi(o) = \bigvee I$ by the preceding paragraph, we conclude
$\psi(x) \neq \psi(y)$.

\textit{Injectivity.}
Let $x, y \in \hat{A}$ with $\phi(x) = \phi(y)$, and suppose $x \neq y$, so
that $x(h) \neq y(h)$ for some $h \in \alg{A}^*$. Since $x$ and $y$ are
locally in $A$, there are $u, v \in A$ with $u(h) = x(h)$ and $v(h) = y(h)$.
Set $p := m(u, x, v)$ and $q := m(u, y, v)$, both in
$[u;v]_{\hat{\alg{A}}}$. Then $p(h) = x(h)$ and $q(h) = y(h)$,
so $p \neq q$, while $\phi(p) = m(u, \phi(x), v) = m(u, \phi(y), v) =
\phi(q)$, contradicting the injectivity of $\phi$ on
$[u;v]_{\hat{\alg{A}}}$. Hence $\phi$ is injective.

\textit{Surjectivity.}
Since $\phi$ fixes $\alg{A}$, the image $\phi(\hat{\alg{A}})$ contains $\alg{A}$.
Since $\hat{\alg{A}}$ is compact and $\alg{A}'$ is Hausdorff, the image
$\phi(\hat{\alg{A}})$ is closed. Since $\alg{A}$ is dense in $\alg{A}'$, the
closed set $\phi(\hat{\alg{A}})$ contains a dense subset of $\alg{A}'$, hence
$\phi(\hat{\alg{A}}) = \alg{A}'$.

Since $\phi$ is a bijective continuous homomorphism from the compact space $\hat{\alg{A}}$
to the Hausdorff space $\alg{A}'$, it is an isomorphism of Boolean topological
median algebras.
\end{proof}

Theorem~\ref{thm:char-profinite} is the internal characterization of
$\hat{\alg{A}}$ announced at the beginning of this section, describing how
$\alg{A}$ sits inside $\hat{\alg{A}}$ through density and compactness,
interval by interval. Globally, the canonical extension of the whole basepoint
order $(\alg{A}, \leq_a)$ also sits inside $\hat{\alg{A}}$, possibly with a
top adjoined, as the set of meets of ideal elements
(Theorem~\ref{thm:delta-inside}).

For a median algebra $\alg{A}$ with finite intervals, the condition of
Theorem~\ref{thm:char-profinite} holds in every $\alg{A}' \in \M_\tau$ in
which $\alg{A}$ is dense. Indeed, the continuity of $z \mapsto m(a, z, b)$
makes $[a;b]_{\alg{A}'}$ the closure of $[a;b]_{\alg{A}}$, as in the first
part of the proof, and this finite bounded distributive lattice is closed and
is its own canonical extension \cite[Corollary~4.31]{FussnerGehrke2025}.
Theorem~\ref{thm:char-profinite} thus gives a second proof of the following
corollary, which is also a consequence of
Proposition~\ref{prop:compactifications}, as we now show.

\begin{corollary}\label{cor:discrete-closure}
Let $\alg{B} \in \M_\tau$ and let $\alg{C}$ be a subalgebra of $\alg{B}$ with
finite intervals. The closure of $\alg{C}$ in $\alg{B}$ is a profinite
completion of $\alg{C}$.
\end{corollary}
\begin{proof}
The closure $\alg{B}'$ of $\alg{C}$ in $\alg{B}$ is a closed subalgebra of
$\alg{B}$, hence lies in $\M_\tau$ (Proposition~\ref{prop:profinite-char}),
and $\alg{C}$ is dense in it. Applied to the inclusion of $\alg{C}$ into
$\alg{B}'$, Proposition~\ref{prop:compactifications} gives an isomorphism
$\alg{B}' \cong P_+$ carrying this inclusion to $e_P$, for a separating
bounded sub-poc set $P$ of $\alg{C}^+$. Since $\alg{C}$ has finite intervals,
the set $P$ is $\alg{C}^+$, and $\alg{B}'$ is isomorphic over $\alg{C}$ to
$(\alg{C}^+)_+$, that is, to $\hat{\alg{C}}$ by
Corollary~\ref{prop:profinite}.
\end{proof}

Thus the profinite completion of a median algebra with finite intervals can be
computed as its closure in any algebra of $\M_\tau$ containing it, for instance
in a power of $\alg{2}_\tau$.

\section{Conditional canonical extensions}\label{sec:ccan}

Without a top, the poset $(\hat{\alg{A}}, \leq_a)$ is usually only a
conditionally complete semilattice (Proposition~\ref{prop:cond-complete}). We
show that it satisfies compactness and one of the two halves of density: every
element is a join of filter elements, but only the elements lying below an
element of $\alg{A}$ are guaranteed to be meets of ideal elements
(Proposition~\ref{prop:global-two-thirds} and Example~\ref{ex:escalator}).

We therefore relax the notion of canonical extension of a poset $\alg{P}$.
Completeness is weakened to conditional completeness, compactness is required
only of the ideals whose join exists, and $\alg{P}$ is asked to be cofinal in
the extension. The resulting conditional canonical extension is unique
(Theorem~\ref{thm:unique-ccan}), and it can be extracted from
$(\hat{\alg{A}}, \leq_a)$ as the union $\hat{\alg{A}}_a$ of the intervals
$[a;b]$ with $b \in A$ (Theorem~\ref{thm:realization}).
Subsection~\ref{subsec:ccan-adjoint} supplies the third point of view on
completions recalled in Section~\ref{sec:can-ext}, a universal property: the
conditional canonical extension is left adjoint to a forgetful functor
(Corollary~\ref{cor:ccan-adjunction}).

\subsection{The conditional canonical extension inside the profinite
completion}\label{subsec:unbounded}

\begin{definition}\label{def:cond-completion}
If $\alg{P}$ is a poset and $\alg{S}$ is a conditionally complete semilattice
such that $\alg{P}$ is a subposet of $\alg{S}$, then $\alg{S}$ is called a
\emph{conditional completion} of $\alg{P}$. 
Filter and ideal elements of
$\alg{S}$, and the notion of density, are defined as in
Definition~\ref{def:can-ext}. 
Any filter $F$ of $\alg{P}$ gives rise to a filter element $\bigwedge F$ in
$\alg{S}$, whereas only ideals $I$ of $\alg{P}$ whose join exists in $\alg{S}$ give
rise to an ideal element $\bigvee I$.
A conditional completion $\alg{S}$ of $\alg{P}$ is \emph{compact} if for every
filter $F$ and every ideal $I$ of $\alg{P}$ such that $\bigvee I$ exists in
$\alg{S}$ and $\bigwedge F \leq \bigvee I$, there are $s\in F$ and $t\in I$ such
that $s \leq t$.
\end{definition}

Conditional completion specializes to the poset completion of \cite{Dunn2005}
exactly when the completing order is a complete lattice.

The profinite completion supplies the first instance of the new notion.

\begin{proposition}\label{prop:cond-complete}
Let $\alg{A}$ be a median algebra.
The semilattice $(\hat{\alg{A}}, \leq_x)$ is conditionally complete for every $x\in
\hat{\alg{A}}$. If in addition $x\in \alg{A}$, then $(\hat{\alg{A}}, \leq_x)$
is a conditional completion of $(\alg{A}, \leq_x)$.
\end{proposition}
\begin{proof}
For the first statement, let $S \subseteq \hat{\alg{A}}$ have an upper bound
$y$ in $(\hat{\alg{A}}, \leq_x)$, and apply Lemma~\ref{lem:ont} to $[x;y]$.
The second statement follows from the first and the fact that
$(\alg{A}, \leq_x)$ is a subposet of $(\hat{\alg{A}}, \leq_x)$.
\end{proof}

The chain $\N$ admits two nonisomorphic
dense and compact conditional completions, namely $\N$ itself and its
canonical extension $\N \cup \{\omega\}$. 
Indeed, the ideal $\N$ has no join
in the former, so that the compactness condition is vacuously satisfied for
this ideal, while its join in the latter is $\omega$. We recover uniqueness by 
the addition of a single order-theoretic condition. 

\begin{definition}\label{def:ccan}
Let $\alg{P}$ be a poset and let $\alg{S}$ be a conditional completion of
$\alg{P}$. The poset $\alg{P}$ is \emph{cofinal} in $\alg{S}$ if every
element of $\alg{S}$ lies below some element of $\alg{P}$. A
\emph{conditional canonical extension} of $\alg{P}$ is a dense and compact
conditional completion of $\alg{P}$ in which $\alg{P}$ is cofinal.
\end{definition}

\begin{theorem}\label{thm:unique-ccan}
A poset has, up to an isomorphism fixing it, at most one conditional
canonical extension.
\end{theorem}
\begin{proof}
 The argument follows the three steps of the proof of
\cite[Theorem~2.5]{Dunn2005}, the cofinality assumption accounting for the one
point at which the present setting departs from the complete-lattice one.
\end{proof}

We now examine how much of Theorem~\ref{thm:can-ext}
survives when the boundedness assumption is dropped, that is,
when the conditional completion $(\hat{\alg{A}}, \leq_a)$ of the whole
basepoint order $(\alg{A}, \leq_a)$ is considered. 

\begin{lemma}\label{lem:dcpo}
Let $\alg{A} \in \M$ and $x \in \hat{A}$. Every up-directed subset $D$ of
$(\hat{\alg{A}}, \leq_x)$ has a join in $(\hat{\alg{A}}, \leq_x)$. This join
is computed coordinatewise in $\alg{2}^{\alg{A}^*}$ and belongs to the
closure of $D$. The dual statement holds for down-directed subsets.
\end{lemma}
\begin{proof}
It is Lemma~\ref{lem:closed-joins}\eqref{it:dir-loc} read in $\hat{\alg{A}}$,
the coordinatewise join of an up-directed subset being locally in that subset
and $\hat{A}$ being closed.
\end{proof}

In particular, every ideal of $(\alg{A}, \leq_a)$ has a join in
$(\hat{\alg{A}}, \leq_a)$, so that the compactness condition is never
satisfied vacuously for this completion.

\begin{proposition}\label{prop:global-two-thirds}
Let $\alg{A} \in \M$ and $a \in A$. The conditional completion
$(\hat{\alg{A}}, \leq_a)$ of $(\alg{A}, \leq_a)$ has the following
properties.
\begin{enumerate}
\item\label{it:gc} It is compact.
\item\label{it:gjd} Every element of $\hat{\alg{A}}$ is a join of filter
elements.
\item\label{it:gmd} Every element of $\hat{\alg{A}}$ that is below some
element of $\alg{A}$ is a meet of ideal elements.
\end{enumerate}
\end{proposition}
\begin{proof}
\eqref{it:gc} By Lemma~\ref{lem:dcpo}, the meet of a filter and the join of
an ideal of $(\alg{A}, \leq_a)$ are computed coordinatewise in
$\alg{2}^{\alg{A}^*}$. The compactness argument in the proof of
Theorem~\ref{thm:can-ext} uses the boundedness of $P$ only through the
coordinatewise computation of these two operations, and therefore applies
with $P$ replaced by $\alg{A}$.

\eqref{it:gjd} Let $x \in \hat{A}$. For every finite subset $F$ of
$\alg{A}^*$, the witness set $G_F = \{c \in A \mid c|_F = x|_F\}$ is
nonempty and down-directed, since $m(a, c, c')$ belongs to $G_F$ whenever
$c$ and $c'$ do. Let $K_F$ be the filter generated by $G_F$ in $(\alg{A},
\leq_a)$. The proof of the identity $x = \bigvee \{\bigwedge K_F \mid F
\Subset \alg{A}^*\}$ in Theorem~\ref{thm:can-ext} is coordinatewise and
uses no joins in $\alg{A}$. The join over the up-directed family
$\{\bigwedge K_F \mid F \Subset \alg{A}^*\}$ exists by
Lemma~\ref{lem:dcpo}. Hence that proof applies with $P$ replaced
by $\alg{A}$.

\eqref{it:gmd} Assume $x \leq_a b$ with $b \in A$. As observed in the proof
of Theorem~\ref{thm:char-profinite}, the interval $[a;b]_{\hat{\alg{A}}}$ is
the closure of $[a;b]_{\alg{A}}$, and is therefore a canonical extension of
$[a;b]_{\alg{A}}$ by Theorem~\ref{thm:can-ext}. By density, the element $x$
is the meet of a family of joins of ideals of $([a;b]_{\alg{A}}, \leq_a)$,
all computed in $[a;b]_{\hat{\alg{A}}}$. Every ideal of $([a;b]_{\alg{A}}, \leq_a)$ is an ideal of $(\alg{A},
\leq_a)$, being an up-directed subset of $[a;b]_{\alg{A}}$ whose
down-closure in $(\alg{A}, \leq_a)$ remains in $[a;b]_{\alg{A}}$. By
Lemma~\ref{lem:ont} its join in $[a;b]_{\hat{\alg{A}}}$ agrees with its join
in $(\hat{\alg{A}}, \leq_a)$. Hence $x$ is a meet of ideal elements of the
completion $(\hat{\alg{A}}, \leq_a)$.
\end{proof}

The following example shows that the restriction in Proposition~\ref{prop:global-two-thirds}\eqref{it:gmd}  cannot be removed in general: the profinite
completion of a median algebra need not be a dense conditional completion of one
of its basepoint orders. We use the fact that, in a compact space, every
sequence has a \emph{cluster point}, that is,  a point every neighbourhood of
which contains infinitely many terms of the sequence.

\begin{example}\label{ex:escalator}
Represent the elements of $\alg{2}^\N \times \alg{2}^\N$ as pairs $(s,t)$
of $\{0,1\}$-sequences, with the median computed coordinatewise.
For natural numbers $n \leq m$, let $c_{n,m}$ be the pair whose first
sequence has $n$ leading ones followed by the constant sequence 0, and whose
second sequence has $m$ leading zeros followed by the constant sequence 1,
\[
  c_{n,m} := \bigl(\, 1 \cdots 1\, 0\, 0 \cdots ,\;
                     0 \cdots 0\, 1\, 1 \cdots \,\bigr).
\]
Let $a$ be the pair of zero sequences, and set
\[
  A = \{a\} \cup \{c_{n,m} \mid n \leq m\}.
\]
A direct computation gives
\begin{gather*}
  m(c_{n_1,m_1}, c_{n_2,m_2}, c_{n_3,m_3})
  = c_{\operatorname{med}(n_1,n_2,n_3),\, \operatorname{med}(m_1,m_2,m_3)},
\\
  m(a, c_{n,m}, c_{n',m'}) = c_{n \wedge n',\, m \vee m'},
\end{gather*}
and the constraint $n \leq m$ is preserved because $\operatorname{med}$ is
monotone.
Thus $A$ is a median subalgebra of $\alg{2}^\N \times \alg{2}^\N$, and $\leq_a$
is the coordinatewise order, that is, $c_{n,m} \leq_a c_{n',m'}$ if and only if
$n \leq n'$ and $m \geq m'$. The diagram below shows the elements $c_{n,m}$ with
$m \leq 4$. The maximal elements $c_{n,n}$ form its top row, and the basepoint
$a$, which lies below all of them, is not drawn.
\begin{center}
\begin{tikzpicture}[x=5mm,y=5mm]
\foreach \n in {0,...,4}
  \foreach \m in {\n,...,4} {
    \fill ({\n+\m},{\n-\m}) circle (1.5pt);
    \ifnum\n<\m \draw ({\n+\m},{\n-\m}) -- ({\n+\m+1},{\n-\m+1}); \fi
    \ifnum\m<4 \draw ({\n+\m},{\n-\m}) -- ({\n+\m+1},{\n-\m-1});
    \else \draw[dotted] ({\n+\m},{\n-\m}) -- ({\n+\m+0.8},{\n-\m-0.8}); \fi
  }
\node[above] at (0,0) {$c_{0,0}$};
\node[above] at (2,0) {$c_{1,1}$};
\node[above] at (4,0) {$c_{2,2}$};
\node[left] at (1,-1) {$c_{0,1}$};
\node[left] at (2,-2) {$c_{0,2}$};
\node[right] at (8.3,0) {$\cdots$};
\end{tikzpicture}
\end{center}

Every ideal of $(\alg{A}, \leq_a)$ is principal. Indeed, let $D$ be an
ideal distinct from $\{a\}$. Whenever $c_{n,m}$ and $c_{n',m'}$ belong to
$D$, the pair has an upper bound $c_{N,M}$ in $D$, and $n \vee n' \leq N \leq
M \leq m \wedge m'$. So, the first subscript of every element of $D$ is
bounded by the second subscript of every element of $D$. The greatest first
subscript $N^*$ and the least second subscript $M^*$ among the elements of
$D$ are therefore attained, and $D$ is the ideal generated by $c_{N^*\!, M^*\!}$. 
In particular, the ideal elements of $(\hat{\alg{A}}, \leq_a)$ are exactly
the elements of $\alg{A}$.

The sequence $(c_{n,n})_{n \in \N}$ has a cluster point $\mu$ in the
compact space $\hat{\alg{A}}$. Fix an index $i$. Once $n > i$, the $i$th
entry of the first sequence of $c_{n,n}$ is $1$, and $0$ in the second.
The restrictions to $A$ of the coordinate projections of
$\alg{2}^\N \times \alg{2}^\N$ belong to $\alg{A}^*$, and evaluation at each of
them is continuous on $\hat{\alg{A}} \subseteq \alg{2}^{\alg{A}^*}$. Hence the
cluster point $\mu$ takes the value $1$ at every coordinate of the first
sequence and $0$ at every coordinate of the second. No element of $\alg{A}$ lies
above $\mu$ in $(\hat{\alg{A}}, \leq_a)$, and in particular $\mu \notin A$. 
So the set of ideal elements above $\mu$ in $(\hat{\alg{A}}, \leq_a)$ is empty,
and $\mu$ is not its meet, which would be a greatest element of
$(\hat{\alg{A}}, \leq_a)$, since $c_{0,0} \not\leq_a \mu$. Hence the
completion $(\hat{\alg{A}}, \leq_a)$ is not dense.

\end{example}

Hence the density of $(\hat{\alg{A}}, \leq_a)$ can fail outside the elements
that lie below an element of $\alg{A}$. The following theorem shows that
restricting the completion to the elements below an element of $\alg{A}$ yields
the conditional canonical extension. 

\begin{theorem}\label{thm:realization}
Let $\alg{A} \in \M$ and $a \in A$, and let
\[
  \hat{\alg{A}}_a \;=\; \bigcup_{b \in A} [a;b]_{\hat{\alg{A}}}
\]
endowed with the subspace topology of
$\hat{\alg{A}}$. 
\begin{enumerate}
\item\label{it:realclosure} 
The set $\hat{\alg{A}}_a$ is a down-set of $(\hat{\alg{A}}, \leq_a)$ and a
topological median subalgebra of $\hat{\alg{A}}$ in which $\alg{A}$ is dense,
and it is dense in $\hat{\alg{A}}$.
\item\label{it:realccan} The poset $(\hat{\alg{A}}_a, \leq_a)$ is the
conditional
canonical extension of $(\alg{A}, \leq_a)$.
\end{enumerate}
\end{theorem}
\begin{proof}
\eqref{it:realclosure} 
The set $\hat{\alg{A}}_a$ is a down-set by construction. The rest of the
statement follows directly from the continuity and monotonicity of $m$, and
the fact that $\alg{A}$ is topologically dense in $\hat{\alg{A}}$.

\eqref{it:realccan} If a subset of $\hat{\alg{A}}_a$ has an upper bound $z$ in
$\hat{\alg{A}}_a$, then $z \leq_a b$ for some $b \in A$ by definition of
$\hat{\alg{A}}_a$, so that the subset is contained in the interval
$[a;b]_{\hat{\alg{A}}}$, which is a complete lattice by Lemma~\ref{lem:ont}. Its
join there is its join in $(\hat{\alg{A}}, \leq_a)$ and again lies in
$\hat{\alg{A}}_a$. Thus $\hat{\alg{A}}_a$ is a conditional completion of
$(\alg{A}, \leq_a)$, in which $\alg{A}$ is cofinal by construction. Every
filter element $\bigwedge F$ of $(\hat{\alg{A}}, \leq_a)$ lies below each
element of $F$, hence is in $\hat{\alg{A}}_a$. For an element $x \leq_a b$ of
$\hat{\alg{A}}_a$, the joins of filters and meets of ideals that
witness density in
Proposition~\ref{prop:global-two-thirds}\eqref{it:gjd}
and~\eqref{it:gmd} are computed inside the interval $[a;b]_{\hat{\alg{A}}}$
and so remain in $\hat{\alg{A}}_a$. Compactness is inherited from
Proposition~\ref{prop:global-two-thirds}\eqref{it:gc}, since the meet of a
filter and the join of a bounded ideal are the same in
$\hat{\alg{A}}_a$ and in $\hat{\alg{A}}$. Hence
$\hat{\alg{A}}_a$ is a conditional canonical extension of
$(\alg{A}, \leq_a)$.
\end{proof}

\begin{example}\label{ex:not-subalgebra}
\emph{The down-set $\hat{\alg{A}}_a$ need not be a subalgebra of
$(\hat{\alg{A}}, \leq_a)$ in the sense of
Definition~\ref{def:cond-complete}.} For instance, let $\alg{A}$ be the
subalgebra $\N^2 \cup \{(\infty, 0), (0, \infty)\}$ of the square of the
chain $\N \cup \{\infty\}$, with basepoint $a = (0,0)$, so that $\leq_a$ is
the product order. The joins $x$ of $\{(n,0) \mid n \in \N\}$ and $y$ of
$\{(0,n) \mid n \in \N\}$ in $(\hat{\alg{A}}, \leq_a)$ exist by
Lemma~\ref{lem:dcpo}, and they lie in $\hat{\alg{A}}_a$, below $(\infty, 0)$
and $(0, \infty)$ respectively. Both lie below the join of $\{(n,n) \mid n
\in \N\}$, so they have a join in $(\hat{\alg{A}}, \leq_a)$. No element of
$\alg{A}$ lies above every $(n,0)$ and every $(0,n)$, so this join is not in
$\hat{\alg{A}}_a$.
\end{example}

By Lemma~\ref{lem:dcpo}, an unbounded up-directed subset of $(\alg{A},
\leq_a)$ still has a join in $(\hat{\alg{A}}, \leq_a)$, and that join lies
outside $\hat{\alg{A}}_a$, since an element of $\alg{A}$ above it would bound
the subset.

By Theorem~\ref{thm:unique-ccan}, the conditional canonical extension
of $(\alg{A}, \leq_a)$ is the whole of $(\hat{\alg{A}}, \leq_a)$ precisely when
$\alg{A}$ is cofinal in $(\hat{\alg{A}}, \leq_a)$,
a condition that holds whenever $(\alg{A}, \leq_a)$ is
upper-bounded, by Corollary~\ref{cor:global-top}.
The following example shows that this sufficient condition is not
necessary.

\begin{example}\label{ex:two-ended-chain}
\emph{The poset $(\alg{A}, \leq_a)$ need not be upper-bounded for $\alg{A}$
to be cofinal in $(\hat{\alg{A}}, \leq_a)$.}
Let $\alg{A}$ be the chain $\Z \cup \{-\infty, +\infty\}$ with its median.
By Corollary~\ref{cor:global-top} for the basepoint $-\infty$, we have
$\hat{\alg{A}} = \alg{A} \cup \{-\infty^+, +\infty^-\}$, where $+\infty^- =
\bigvee \N$ and $-\infty^+ = \bigwedge (-\N)$ in the order of the chain. For the
basepoint $0$ the order $\leq_0$ on $\hat{\alg{A}}$ consists of two chains 
which meet only at $0$, and $(\alg{A}, \leq_0)$ is not upper-bounded. Every element
of $\hat{\alg{A}}$ nevertheless lies below $+\infty$ or below $-\infty$, so
that $\alg{A}$ is cofinal in $(\hat{\alg{A}}, \leq_0)$, and $(\hat{\alg{A}},
\leq_0)$ is the conditional canonical extension of $(\alg{A}, \leq_0)$.
\end{example}

The canonical extension of the whole basepoint order, in the sense of
Definition~\ref{def:can-ext}, is a complete lattice, unlike the conditional
one, and it can also be located inside $\hat{\alg{A}}$. It
rests on the following observation on posets.

\begin{lemma}\label{lem:ce-inside}
Let $\alg{P}$ be a poset and let $\alg{L}$ be a compact completion of
$\alg{P}$ in which every element is a join of filter elements. The set $M$ of
the elements of $\alg{L}$ that are meets of ideal elements, as a subposet of
$\alg{L}$, is a canonical extension of $\alg{P}$.
\end{lemma}
\begin{proof}
The proof is that of \cite[Theorem~2.8]{Gouveia2013}, where $\alg{L}$ is the iterated
filter lattice of a unital meet-semilattice. It applies verbatim, since it uses
only that $M$ is closed under the meets of $\alg{L}$, that the filter and ideal
elements of $M$ are those of $\alg{L}$, and that $\alg{P}$ lies in $M$, each
element of $\alg{P}$ being an ideal element.
\end{proof}

\begin{theorem}\label{thm:delta-inside}
Let $\alg{A} \in \M$ and $a \in A$. Let $\alg{L}$ be $(\hat{\alg{A}},
\leq_a)$ if this poset has a greatest element, and the poset obtained from
$(\hat{\alg{A}}, \leq_a)$ by adding a fresh greatest element $\top$
otherwise.
\begin{enumerate}
\item\label{it:di-top} The poset $(\hat{\alg{A}}, \leq_a)$ has a greatest
element if and only if $\alg{A}$ has an upper bound in it, in particular if
$(\alg{A}, \leq_a)$ is up-directed.
\item\label{it:di-where} The canonical extension of $(\alg{A}, \leq_a)$ is,
up to an isomorphism fixing $A$, the set of the elements of $\alg{L}$ that are
meets of ideal elements.
\item\label{it:di-between} This set contains $\hat{\alg{A}}_a$, and it is the
whole of $\alg{L}$ if and only if every element of $\hat{\alg{A}}$ is a meet
of ideal elements of $(\hat{\alg{A}}, \leq_a)$.
\end{enumerate}
\end{theorem}
\begin{proof}
\eqref{it:di-top} Let $z \in \hat{A}$ be an upper bound of $A$. Every element
of $\hat{\alg{A}}$ is a join of filter elements by
Proposition~\ref{prop:global-two-thirds}\eqref{it:gjd}, each of them below
$z$, so $z$ is the greatest element. The converse is clear, and if $(\alg{A},
\leq_a)$ is up-directed, then its join exists by Lemma~\ref{lem:dcpo}.

\eqref{it:di-where} By Lemma~\ref{lem:ce-inside}, it suffices to show that
$\alg{L}$ is a compact completion of $(\alg{A}, \leq_a)$ in which every element
is a join of filter elements. It is a complete lattice by
Proposition~\ref{prop:cond-complete}, and its filter and ideal elements are
those of $(\hat{\alg{A}}, \leq_a)$, so both properties follow from
Proposition~\ref{prop:global-two-thirds}, the element $\top$, when present,
being the join of $A$ by~\eqref{it:di-top}.

\eqref{it:di-between} Every element of $\hat{\alg{A}}_a$ lies below an
element of $\alg{A}$, so it is a meet of ideal elements by
Proposition~\ref{prop:global-two-thirds}\eqref{it:gmd}. The second statement
follows, since the greatest element of $\alg{L}$ is the empty meet and the
meets of nonempty sets of ideal elements are computed in $(\hat{\alg{A}},
\leq_a)$.
\end{proof}

When $(\alg{A}, \leq_a)$ has a greatest element, every element of
$\hat{\alg{A}}$ lies below an element of $\alg{A}$, and
Theorem~\ref{thm:delta-inside} returns Corollary~\ref{cor:global-top}. In
general the canonical extension of $(\alg{A}, \leq_a)$ lies between
$\hat{\alg{A}}_a$ and $\alg{L}$, and the second inclusion can be strict, as the
element $\mu$ of Example~\ref{ex:escalator} is not a meet of ideal elements.
For a unital meet-semilattice $S$, Gouveia and Priestley identify the canonical
extension $S^\delta$ of $S$ in
the same way with the meets of directed joins inside the profinite completion
of $S$ in the variety of semilattices \cite[Theorem~3.2]{Gouveia2013}.
Theorem~\ref{thm:delta-inside} does the same inside the profinite completion
in $\M$, where the basepoint order need not have a top.

\begin{remark}\label{rem:ccan-vs-zero}
The conditional canonical extension and the zero-completion recalled in
Subsection~\ref{subsec:roller} are both extensions of $\alg{A}$ built from
its intervals, and both sit inside $\hat{\alg{A}}$. The zero-completion is
the inverse limit of the intervals $[u;v]_{\alg{A}}$ under the gate
projections $x \mapsto m(u,v,x)$ \cite[\S4.1]{Fioravanti2020}, and it embeds in
$\hat{\alg{A}}$ over $\alg{A}$ \cite[Lemma~4.3(2)]{Fioravanti2020}. It adds
points outside the intervals of $\alg{A}$ and leaves these intervals unchanged,
since $\alg{A}$ is convex in it \cite[Theorem~1(i)]{BandeltMeletiou1993},
\cite[Corollary~4.4(1)]{Fioravanti2020}. The conditional canonical extension
$\hat{\alg{A}}_a$ is the colimit of the canonical extensions
$[a;b]_{\hat{\alg{A}}}$ of the intervals $[a;b]_{\alg{A}}$ with $b \in A$
(Theorem~\ref{thm:char-profinite}) under the inclusions
$[a;b]_{\hat{\alg{A}}} \subseteq [a;b']_{\hat{\alg{A}}}$ for $b \leq_a b'$.
It does the opposite, adding no point outside these intervals.
Consequently the zero-completion meets $\hat{\alg{A}}_a$ exactly in
$\alg{A}$, because an element of both lies in an interval $[a;b]$ of the
zero-completion with $b \in A$, and that interval is $[a;b]_{\alg{A}}$.

\end{remark}

\subsection{The conditional canonical extension as a left adjoint}
\label{subsec:ccan-adjoint}

We now show that the conditional canonical extension is a left adjoint.

\begin{definition}\label{def:cat-ms}
Let $\cate{MS}$ be the category whose objects are the median semilattices and
whose morphisms are the maps that preserve the bottom element, the binary meets,
and all existing binary joins. Let $\cate{M}_\bullet$ be the category of pointed
median algebras $(\alg{A}, a)$ with the median homomorphisms that preserve the
basepoint.
\end{definition}

Sholander's correspondence, recalled after
Definition~\ref{def:median-semilattice}, puts the objects of the two
categories in bijection. It also matches their morphisms, so that the two
categories are isomorphic over the underlying sets.

\begin{proposition}\label{prop:ms-pointed}
The assignments $(\alg{A}, a) \mapsto (\alg{A}, \wedge_a)$ and
$(S, \wedge) \mapsto (\alg{S}, 0)$ are mutually inverse
isomorphisms between $\cate{M}_\bullet$ and $\cate{MS}$, and they act as the
identity on morphisms.
\end{proposition}
\begin{proof}
The meet $\wedge_a$ is the term $m(a, -, -)$, and by
Lemma~\ref{lem:interval-dual} the join of two elements $x, y$ below some $c$
is the term $m(x, y, c)$, so basepoint-preserving median homomorphisms are
morphisms of $\cate{MS}$. Conversely, by \eqref{eq:semi-to-med} the median is
built from binary meets and joins of pairs bounded above, which morphisms of
$\cate{MS}$ preserve.
\end{proof}

\begin{definition}\label{def:cat-cms}
A median semilattice is \emph{conditionally completely distributive} if it is
conditionally complete and each of its principal ideals is a completely
distributive lattice. Let $\cate{CMS}$ be the category of conditionally
completely distributive median semilattices, whose morphisms are the maps
preserving all nonempty meets and all existing joins, the empty join included.
Write $U_\Delta \colon \cate{CMS} \to \cate{MS}$ for the forgetful functor.
\end{definition}

The next lemma shows that the objects of the preceding sections belong to
$\cate{CMS}$.

\begin{lemma}\label{lem:mtau-cms}
Let $\alg{X}, \alg{Y} \in \M_\tau$,  $x \in X$ and $y\in Y$.
\begin{enumerate}
\item\label{it:mc-obj} The poset $(\alg{X}, \leq_x)$ belongs to 
$\cate{CMS}$.

\item\label{it:mc-mor} If $f \in \cate{M}_\tau(\alg{X}, \alg{Y})$
is such that $f(x) = y$, then $f$ is a $\cate{CMS}$-morphism
from $(\alg{X}, \leq_x)$ to $(\alg{Y}, \leq_y)$.

\item\label{it:mc-hat} For $\alg{A} \in \M$ and $a \in A$, both
$(\hat{\alg{A}}, \leq_a)$ and $\hat{\alg{A}}_a$ belong to $\cate{CMS}$,
and the inclusion of $(\alg{A}, \leq_a)$ into either of them is a
$\cate{MS}$-morphism.
\end{enumerate}
\end{lemma}
\begin{proof}
\eqref{it:mc-obj} The proof is an application of Lemma~\ref{lem:ont}.

\eqref{it:mc-mor} By Proposition~\ref{prop:ms-pointed}, the map $f$ preserves
the bottom, the binary meets and the existing binary joins. Let $D$ be an
up-directed subset of $\alg{X}$ bounded above, with join $s$, and let $t$ be
the join of $f(D)$ in $(\alg{Y}, \leq_y)$. Monotonicity gives $t \leq_y
f(s)$. The principal ideal $[y;t]$ is closed in $\alg{Y}$, so $f^{-1}([y;t])$
is a closed subset of $\alg{X}$ containing $D$. Since $s$ lies in the closure of
$D$ by Lemma~\ref{lem:closed-joins}\eqref{it:dir-loc}, it also lies in
$f^{-1}([y;t])$, which shows that $f(\bigvee D) = \bigvee f(D)$.
Since a bounded subset has the same join as the up-directed set of the joins
of its nonempty finite subsets, the map $f$ preserves all existing joins. The
argument for nonempty meets is dual, with the principal filter ${\uparrow}t =
\{w \mid m(y,t,w) = t\}$, which is also closed, in place of $[y;t]$.

\eqref{it:mc-hat} Apply \eqref{it:mc-obj} to $\alg{X} = \hat{\alg{A}}$. The
poset $\hat{\alg{A}}_a$ is conditionally complete by
Theorem~\ref{thm:realization}\eqref{it:realccan}, and its principal ideals are
the intervals $[a;b]_{\hat{\alg{A}}}$ with $b \in \hat{\alg{A}}_a$, and are
completely
distributive by Lemma~\ref{lem:ont}. Both inclusions of $(\alg{A}, \leq_a)$
are basepoint-preserving median homomorphisms, hence morphisms of $\cate{MS}$
by Proposition~\ref{prop:ms-pointed}.
\end{proof}

Writing $\alg{S}^\Delta$ for the conditional canonical extension of
$\alg{S} \in \cate{MS}$, that is, $\hat{\alg{A}}_a$ when $\alg{S}
= (\alg{A}, \leq_a)$, we can state the universal property that we are looking
for.

\begin{theorem}\label{thm:ccan-free}
Let $\alg{S} \in \cate{MS}$ and $\alg{C} \in \cate{CMS}$. Every morphism $h
\colon \alg{S} \to U_\Delta(\alg{C})$ of $\cate{MS}$ extends to a unique
morphism $\tilde{h} \colon \alg{S}^\Delta \to \alg{C}$ of $\cate{CMS}$.
\end{theorem}
\begin{proof}
Write $\alg{S} = (\alg{A}, \leq_a)$, so that $\alg{S}^\Delta =
\hat{\alg{A}}_a = \bigcup_{b \in A} [a;b]_{\hat{\alg{A}}}$. We first extend
$h$ to each of these intervals, and then glue the extensions into a single
map on $\hat{\alg{A}}_a$.

\textit{The local extensions.} Fix $b \in A$. The interval $[a;b]_{\alg{A}}$
is the principal ideal ${\downarrow}b$ of $\alg{S}$, so it is a bounded
distributive lattice, and $h$ restricts to a homomorphism of bounded lattices
from $[a;b]_{\alg{A}}$ to $[0;h(b)]_{\alg{C}}$, where $0$ is the bottom of
$\alg{C}$. The codomain is a complete lattice, since it is a principal ideal
of a conditionally complete semilattice, and it is completely distributive by
hypothesis. By Theorem~\ref{thm:char-profinite} the interval
$[a;b]_{\hat{\alg{A}}}$ is the canonical extension of $[a;b]_{\alg{A}}$, and
the canonical extension of a bounded distributive lattice is free over it
among the completely distributive lattices \cite{MortonvanAlten2018} (see
also \cite[Theorem~2.1]{BezhanishviliHardingJibladze2021};
\cite[Theorem~8.21]{FussnerGehrke2025}). There is therefore a unique
complete lattice homomorphism
\[
  h_b \colon [a;b]_{\hat{\alg{A}}} \longrightarrow [0;h(b)]_{\alg{C}}
\]
extending the restriction of $h$.

\textit{Coherence.} Let $b, b' \in A$. We show that $h_b=h_{b'}$ on
$[a;b]_{\hat{\alg{A}}} \cap
[a;b']_{\hat{\alg{A}}}$, which is $[a;c]_{\hat{\alg{A}}}$ for $c := m(a, b,
b')$. By Lemma~\ref{lem:ont} these
intervals are complete subalgebras of $(\hat{\alg{A}}, \leq_a)$, so the
restriction $g$ of $h_b$ to $[a;c]_{\hat{\alg{A}}}$ preserves all nonempty
meets and joins. Since $h_b$ is order preserving and $h_b(c) = h(c)$, the map
$g$ takes its values in $[0;h(c)]_{\alg{C}}$, whose nonempty meets and joins
are those of $\alg{C}$. Hence $g$ is a homomorphism of bounded complete
lattices extending the restriction of $h$. It equals $h_c$ by the uniqueness
mentioned in the first part of the proof, and likewise for $h_{b'}$. Setting
$\tilde{h}(z) := h_b(z)$ for any $b \in A$ with $z \leq_a b$ therefore defines
a map on $\hat{\alg{A}}_a$ extending $h$.

\textit{Binary meets.} Let $z \leq_a b$ and $z' \leq_a b'$ with $b, b' \in A$.
The element $z \wedge z'$ equals $z \wedge (z' \wedge b)$, a meet of two
elements of $[a;b]_{\hat{\alg{A}}}$, and $z' \wedge b$ equals $z' \wedge m(a,
b, b')$, a meet of two elements of $[a;b']_{\hat{\alg{A}}}$. Applying $h_b$
and then $h_{b'}$ gives
\[
  \tilde{h}(z \wedge z') = \tilde{h}(z) \wedge \tilde{h}(z') \wedge h(b)
  \wedge h(b') = \tilde{h}(z) \wedge \tilde{h}(z'),
\]
since $\tilde{h}(z) \leq h(b)$ and
$\tilde{h}(z') \leq h(b')$. Hence $\tilde{h}$ preserves binary meets, and it
is in particular order preserving.

\textit{Arbitrary meets and joins.} Let $S$ be a nonempty subset of
$\hat{\alg{A}}_a$, choose $s_0 \in S$, and choose $b \in A$ with
$s_0 \leq_a b$. The meet of $S$ is the meet of the subset $\{s \wedge s_0 \mid
s \in S\}$ of $[a;b]_{\hat{\alg{A}}}$, so that
\[
  \tilde{h}\Bigl(\bigwedge S\Bigr)
  = \bigwedge_{s \in S} \bigl(\tilde{h}(s) \wedge \tilde{h}(s_0)\bigr)
  = \bigwedge \tilde{h}(S),
\]
the first equality by the preceding paragraph and the fact that $h_b$
preserves meets, the second because $\tilde{h}(s_0)$ is one of the terms. A
subset of $\hat{\alg{A}}_a$ with an upper bound has one in
$\alg{A}$, hence lies in a single interval $[a;b]_{\hat{\alg{A}}}$, where its
join is computed by Lemma~\ref{lem:ont} and preserved by $h_b$. The empty join
is the bottom, and so is preserved by $\tilde{h}$.

\textit{Uniqueness.} Every element of $\hat{\alg{A}}_a$ lies in an
interval $[a;b]_{\hat{\alg{A}}}$ with $b \in A$, where it is a join of meets
of elements of $[a;b]_{\alg{A}}$, and these meets and joins are those of
$\hat{\alg{A}}_a$ by Lemma~\ref{lem:ont}. A morphism of
$\cate{CMS}$ agreeing with $h$ on $\alg{A}$ therefore agrees with $\tilde{h}$.
\end{proof}

Theorem~\ref{thm:ccan-free} makes the conditional canonical extension a
functor. A morphism $f \colon \alg{S} \to \alg{T}$ of $\cate{MS}$ is sent to
the unique morphism $f^\Delta \colon \alg{S}^\Delta \to \alg{T}^\Delta$ of
$\cate{CMS}$ extending $f$, and the uniqueness gives $\mathrm{id}^\Delta =
\mathrm{id}$ and $(g \circ f)^\Delta = g^\Delta \circ f^\Delta$.

\begin{corollary}[The conditional canonical extension adjunction]
\label{cor:ccan-adjunction}
The conditional canonical extension functor $(-)^\Delta$ is left adjoint to
$U_\Delta$. In other words, for every $\alg{S} \in \cate{MS}$ and $\alg{C}
\in \cate{CMS}$ there is a bijection
\[
  \cate{MS}\bigl(\alg{S}, U_\Delta(\alg{C})\bigr) \;\cong\;
  \cate{CMS}\bigl(\alg{S}^\Delta, \alg{C}\bigr),
\]
natural in $\alg{S}$ and $\alg{C}$. The unit of the adjunction is the
inclusion of $\alg{S}$ into $\alg{S}^\Delta$.
\end{corollary}

The following result characterizes the median semilattices that are their own
conditional canonical extension, and shows that this property does not depend
on the basepoint.

\begin{proposition}\label{prop:ccan-unit-iso}
Let $\alg{S} = (\alg{A}, \leq_a) \in \cate{MS}$. The unit of
Corollary~\ref{cor:ccan-adjunction} at $\alg{S}$ is an isomorphism if and only
if $\alg{A}$ is discrete, that is, has finite intervals.
\end{proposition}
\begin{proof}
The unit is an isomorphism if and only if $[a;b]_{\hat{\alg{A}}} =
[a;b]_{\alg{A}}$ for every $b \in A$, that is, by
Theorem~\ref{thm:char-profinite}, if and only if every bounded distributive
lattice $[a;b]_{\alg{A}}$ is its own canonical extension. This happens if and
only if each of them is finite \cite[Corollary~4.31]{FussnerGehrke2025}, that
is, by Lemma~\ref{lem:discrete-ideals}, if and only if $\alg{A}$ is discrete.
\end{proof}

\section{Transferring properties to the Roller compactification}
\label{sec:transfer}

Several results of the geometric literature bound the rank of the Roller
compactification $\hat{\alg{A}}$ of a median algebra $\alg{A}$, or of the
components of its boundary. Their non-strict inequalities all follow from
\cite[Lemmas~8.2.1 and~9.3.1]{Bowditch2022}, since these components are
subalgebras of $\hat{\alg{A}}$ and rank does not increase in subalgebras.
This section gives a logical reading of these preservation results, as an
application of a general transfer theorem (Theorem~\ref{thm:transfer}).

For $n \in \N$, an \emph{$n$-cube} of a median algebra $\alg{A}$ is a subalgebra
isomorphic to $\alg{2}^n$. The \emph{rank} of $\alg{A}$, denoted by
$\rank(\alg{A})$, is the largest $n$ for which $\alg{A}$ has an $n$-cube, and
$\rank(\alg{A}) = \infty$ if there is no largest one \cite[\S8.2]{Bowditch2022}.
A median algebra has rank $0$ exactly when it has one element. Rank does not
increase in subalgebras or homomorphic images \cite[\S8.2]{Bowditch2022}. The
first fact is immediate, and Lemma~\ref{lem:cube-lift} below derives the second
from the duality of Theorem~\ref{thm:dualities}.

\begin{lemma}\label{lem:cube-lift}
For every $n \in \N$, the dual $(\alg{2}^n)^*$ is injective in $\Poc_\tau$.
Consequently, if $f \colon \alg{A} \twoheadrightarrow \alg{B}$ is a surjective
homomorphism, then $\rank(\alg{B}) \leq \rank(\alg{A})$.
\end{lemma}
\begin{proof}
Let $S$ be the underlying set of $\alg{2}^n$, let $\mathbf{0}$ be the constant
tuple $0$ of $S$, and let $\alg{P} := (\alg{2}^n)^*$, seen as a closed
substructure of $(\utilde{\alg{2}}_\tau)^S$. The set $Q := P \setminus \{0,1\}$
consists of the $n$ projections and their complements, and it is an antichain
closed under ${}^\co$. For a map $w \colon S \to \{0,1\}$ that is not in $Q$,
let $c(w)$ be $1$ if $p < w$ for some
$p \in Q$, be $0$ if $w < p$ for some $p \in Q$, and be $w(\mathbf{0})$
otherwise. No $w$ lies strictly between two members of the antichain $Q$, so $c$
is well defined. Let $r \colon (\utilde{\alg{2}}_\tau)^S \to \alg{P}$ be the map
that fixes $Q$ and sends every $w \notin Q$ to the constant map with value
$c(w)$. It is continuous, fixes $P$ and commutes with ${}^\co$. It also
preserves the order. Let $w < w'$. If $w \in Q$, then $w' \notin Q$, as $Q$ is
an antichain, and $r(w') = 1$, and symmetrically $r(w) = 0$ if $w' \in Q$.
Otherwise, it suffices to show that $c(w) = 1$ implies $c(w') = 1$. Either $w$
lies above a member of $Q$, and so does $w'$, or $w(\mathbf{0}) = 1$ and $w$
lies below no member of $Q$. In the second case $w'$ lies below no member of
$Q$ either, and $w'(\mathbf{0}) \geq w(\mathbf{0}) = 1$, so that $c(w') = 1$.
Hence $\alg{P}$ is a retract of
$(\utilde{\alg{2}}_\tau)^S$, which is injective in $\Poc_\tau$ as a power of
$\utilde{\alg{2}}_\tau$, since the duality of
Theorem~\ref{thm:dualities}(\ref{it:dual-poctau}) is strong. So $\alg{P}$ is
injective.

For the second assertion, let $u \colon \alg{2}^n \to \alg{B}$ be an injective
homomorphism. Since $f$ is surjective, the morphism $f^*$ is an embedding, and
since $\alg{P}$ is injective and the duality is full, there is a homomorphism
$v \colon \alg{2}^n \to \alg{A}$ with $v^* \circ f^* = u^*$. Then $f \circ v =
u$ by faithfulness, so $v$ is injective, and $f$ maps the $n$-cube
$v(\alg{2}^n)$ of $\alg{A}$ onto $u(\alg{2}^n)$.
\end{proof}

Let $L_{\M}$ be the first-order language with a single ternary function symbol
$m$, so that median algebras are $L_{\M}$-structures. A sentence $\varphi$ of
$L_{\M}$ is \emph{universal} if it has the form $\forall x_1 \cdots \forall
x_n\, \psi$ with $\psi$ quantifier-free. A universal sentence that holds in an
algebra trivially holds in each of its subalgebras. If in addition the
formula $\psi$ is
\emph{positive}, that is, built from atomic formulas by conjunctions and
disjunctions, then $\varphi$ is called \emph{positive universal}. Positive
universal sentences are preserved by homomorphic images, since positive
quantifier-free formulas are preserved by homomorphisms. We denote
by $\Mod(\varphi)$ the class of median algebras that satisfy a sentence
$\varphi$ of $L_{\M}$, and for a class $\clK$ of median algebras we write
$\clK \models \varphi$ if every member of $\clK$ satisfies $\varphi$. The two
facts above thus read
$\classop{S}(\Mod(\varphi)) = \Mod(\varphi)$ for a universal sentence
$\varphi$, and $\classop{H}(\Mod(\varphi)) = \Mod(\varphi)$ if $\varphi$ is
moreover positive.
A first example is the positive universal sentence
\[
  \lambda := \forall x \forall y \forall z\, \bigl(m(x,y,z) = x \vee
  m(x,y,z) = y \vee m(x,y,z) = z\bigr),
\]
which defines the \emph{conservative} median algebras, namely the median
algebras of totally ordered sets and the square $\alg{2}^2$
\cite[Theorem~3]{CouceiroMarichalTeheux2016}.

\begin{example}\label{ex:rho-k}
For a finite median algebra $\alg{B}$, with a fresh variable $x_b$ for every
$b \in B$, let
\[
  \varphi_{\alg{B}} := \forall (x_b)_{b \in B} \Bigl( \bigwedge_{b \neq c}
  x_b \neq x_c \;\to\; \bigvee_{b, c, d \in B} m(x_b, x_c, x_d) \neq
  x_{m(b,c,d)} \Bigr).
\]
An assignment in a median algebra $\alg{A}$ falsifies the quantifier-free part
of $\varphi_{\alg{B}}$ exactly when $b \mapsto x_b$ is an injective
homomorphism from $\alg{B}$ to $\alg{A}$, so that $\alg{A}$ satisfies
$\varphi_{\alg{B}}$ if and only if none of its subalgebras is isomorphic to
$\alg{B}$. For $k \geq 0$, the sentence $\rho_k := \varphi_{\alg{2}^{k+1}}$
thus satisfies, for every median algebra $\alg{A}$,
\[
\alg{A} \models \rho_{k} \quad \iff \quad \rank(\alg{A}) \leq k.
\]
Although $\rho_k$ is not positive universal as written,
Lemma~\ref{lem:cube-lift} gives the equality
$\classop{H}(\Mod(\rho_k)) = \Mod(\rho_k)$. For $k \geq 1$,
the class $\Mod(\rho_k)$ is not closed under products, since
$\alg{2}^{k+1}$ is a product of $k+1$ copies of $\alg{2}$.
\end{example}

The finite quotients of $\alg{A}$ are its quotients $\alg{A}/\theta_F$ by a
finite set $F$ of prime convex subsets, that is, the finite
hyperplane-collapses when $\alg{A}$ is the vertex set of a median graph
(Subsection~\ref{subsec:roller}). We show that they capture the finitely
generated subalgebras of $\hat{\alg{A}}$.

\begin{lemma}\label{lem:finite-sub}
Let $\alg{A}$ be a median algebra. Every finitely generated subalgebra of
$\hat{\alg{A}}$ embeds into a finite quotient of $\alg{A}$ by a limit
projection.
\end{lemma}
\begin{proof}
Let $\alg{B}$ be a finitely generated subalgebra of $\hat{\alg{A}}$. It is
finite, the variety $\M$ being locally finite. For each pair
of distinct elements of $B$, choose a congruence of finite index at whose
coordinate the two threads differ. The intersection $\alpha$ of these finitely
many congruences has finite index. The limit projection $p_\alpha \colon
\hat{\alg{A}} \to \alg{A}/\alpha$ is injective on $B$, and its restriction to
$\alg{B}$ is an embedding.
\end{proof}

The next theorem combines basic facts of model theory with
Lemma~\ref{lem:finite-sub} to transfer universal sentences from the finite
quotients of a median algebra to its profinite completion.

\begin{theorem}\label{thm:transfer}
Let $\alg{A}$ be a median algebra and let $\varphi$ be a universal sentence
of $L_{\M}$.
\begin{enumerate}
\item\label{it:transfer-quot} If $\varphi$ holds in every finite quotient of
  $\alg{A}$, then $\classop{S}(\hat{\alg{A}}) \models \varphi$.
\item\label{it:transfer-closed} If $\alg{A} \models \varphi$ and
  $\classop{H}(\Mod(\varphi)) = \Mod(\varphi)$, then
  $\classop{SH}(\hat{\alg{A}}) \models \varphi$.
\item\label{it:transfer-pos} In particular, if $\varphi$ is positive universal
  and $\alg{A} \models \varphi$, then $\classop{SH}(\hat{\alg{A}}) \models
  \varphi$.
\end{enumerate}
\end{theorem}
\begin{proof}
For~(\ref{it:transfer-quot}), it suffices to show that $\varphi$ holds in
$\hat{\alg{A}}$.
Suppose that $\varphi$ fails in $\hat{\alg{A}}$ at a tuple $\mathbf{x}$. It then
fails in the subalgebra generated by $\mathbf{x}$, which is isomorphic, by
Lemma~\ref{lem:finite-sub}, to a subalgebra of a finite quotient of $\alg{A}$,
and so in that finite quotient.
For~(\ref{it:transfer-closed}), the closure assumption makes $\varphi$ hold in
every finite quotient of $\alg{A}$. Hence it holds in $\hat{\alg{A}}$
by~(\ref{it:transfer-quot}), and in (every subalgebra of) every homomorphic
image of $\hat{\alg{A}}$. Item~(\ref{it:transfer-pos}) is the instance
of~(\ref{it:transfer-closed}) for positive universal sentences.
\end{proof}

\begin{remark}\label{rem:transfer-not-necessary}
The hypothesis of Theorem~\ref{thm:transfer}(\ref{it:transfer-quot}) is
sufficient but not necessary, because a finite quotient of $\alg{A}$ need not be
isomorphic to a subalgebra of $\alg{A}$. Consider the universal sentence
\[
  \sigma := \forall x_0 \cdots \forall x_4 \Bigl(
  \bigwedge_{0 \leq i < j \leq 4} x_i \neq x_j \;\to\;
  \bigvee_{1 \leq i < j < k \leq 4} m(x_i, x_j, x_k) \neq x_0 \Bigr),
\]
which forbids five distinct elements $x_0, \ldots, x_4$ such that $x_0$ is the
median of any three of $x_1, \ldots, x_4$. It holds in the vertex set
$\alg{D}$ of the tree obtained from an edge $uv$ by attaching two leaves to
each of $u$ and $v$, since $x_1, \ldots, x_4$ would lie in four distinct
branches at $x_0$, whereas no vertex of $\alg{D}$ has more than three. It
fails in the quotient of $\alg{D}$ that identifies $u$ with $v$ and nothing
else, which is the star with four leaves, and which is the finite quotient
$\alg{D}/\theta_F$ for the set $F$ of the prime convex subsets of $\alg{D}$
that do not separate $u$ from $v$. As $\alg{D}$ is finite, it is its own
profinite completion, so that $\sigma$ holds in $\hat{\alg{D}}$ although it
fails in a finite quotient of $\alg{D}$.
\begin{center}
\begin{tikzpicture}[x=6mm,y=6mm,every node/.style={font=\small}]
\foreach \p in {(0,0),(2,0),(-1,1),(-1,-1),(3,1),(3,-1)}
  \fill \p circle (1.5pt);
\draw (-1,1) -- (0,0) -- (-1,-1) (0,0) -- (2,0) (3,1) -- (2,0) -- (3,-1);
\node[below] at (0,0) {$u$}; \node[below] at (2,0) {$v$};
\node at (1,-2.2) {$\alg{D}$};
\draw[->] (4,0) -- (6,0);
\foreach \p in {(8,0),(7,1),(7,-1),(9,1),(9,-1)} \fill \p circle (1.5pt);
\draw (7,1) -- (9,-1) (7,-1) -- (9,1);
\node at (8,-2.2) {$\alg{D}/\theta_F$};
\end{tikzpicture}
\end{center}
\end{remark}

Note that the members of $\classop{SH}(\hat{\alg{A}})$ that appear in
Theorem~\ref{thm:transfer} include several of the extensions of
$\alg{A}$ met in the geometric literature.  Applying
Theorem~\ref{thm:transfer}(\ref{it:transfer-closed}) to the sentences $\rho_k$
of Example~\ref{ex:rho-k}, we obtain the following corollary.

\begin{corollary}\label{cor:rank-transfer}
Let $\alg{A}$ be a median algebra.
\begin{enumerate}
\item Every $\alg{B} \in \classop{SH}(\hat{\alg{A}})$ satisfies $\rank(\alg{B})
  \leq \rank(\alg{A})$.
\item The profinite completion $\hat{\alg{A}}$, the zero-completion of
  $\alg{A}$, the conditional canonical extensions $\hat{\alg{A}}_a$ for $a \in
  A$, and every Boolean topological median algebra in which $\alg{A}$ is dense
  have the same rank as $\alg{A}$.
\item Every component of the Roller compactification of a median graph, or of a
  complete median space of finite rank, has rank at most that of the graph or of
  the space.
\end{enumerate}
\end{corollary}

Item~(2) is \cite[Lemma~12.3.4]{Bowditch2022}, which covers every
topological median algebra with a dense subalgebra, together with
\cite[Lemma~4.9]{Fioravanti2020} for the zero-completion, and item~(3) is the
non-strict half of \cite[Proposition~2.2.2]{Genevois2024} and of
\cite[Theorem~D]{Fioravanti2020}.

For instance, the profinite completion of the median algebra $\alg{T}$ of the
infinite rooted binary tree, whose vertices are the words of
$\{0,1\}^{<\omega}$, is the set $\{0,1\}^{<\omega} \cup \{0,1\}^\omega$ of
the finite and infinite words, by Roller's description of the double dual of a
tree \cite[Examples~4.3(vi)]{Roller2016} and Corollary~\ref{prop:profinite}.
By Corollary~\ref{cor:rank-transfer}, it has rank one, although its ends do
not form a subalgebra, the median of three ends being the vertex at which they
part.

\section{Conclusion and perspectives}\label{sec:conclusion}

Universal algebra has delivered the dualities, the profinite completion and
its order-theoretic descriptions together, and has presented the completions
as universal constructions. We hope that these results will prove to be more
than an elegant repackaging, and that they will become useful tools for the
geometric group theory community. The transfer principle of
Section~\ref{sec:transfer} is a first instance, in which the description of
the Roller compactification as a limit of finite quotients turns rank bounds
of the geometric literature into consequences of a single preservation
theorem. In the tradition of universal algebra and algebraic logic, we close
with several routes for further research.

\bigskip

\paragraph{\bfseries Correspondence theory}
In modal logic, correspondence theory relates the validity of a formula on a
frame, a monadic second-order condition obtained through the standard
translation, to first-order properties of the frame
\cite{BlackburnDeRijkeVenema2001}. The construction of the
complex algebra makes bounded poc sets a natural semantics for $L_{\M}$, and a
similar translation can be considered.

Let $\alg{X}$ be a bounded poc set. A \emph{valuation} on $\alg{X}$ is a map
$v$ that assigns to each variable $x$ of $L_{\M}$ a morphism $v(x) \colon
\alg{X} \to \utilde{\alg{2}}$. The \emph{value} of a term $t$ under $v$ is the
map $t^v \colon X \to \{0, 1\}$ defined by $x^v := v(x)$ and
\[
  m(t_1, t_2, t_3)^v(p) := m\bigl(t_1^v(p), t_2^v(p), t_3^v(p)\bigr)
  \qquad (p \in X),
\]
where $m$ on the right is the majority operation of $\{0, 1\}$. As for the
complex algebra (Definition~\ref{def:complex}), each $t^v$ is again a
morphism. We write
$(\alg{X}, v) \Vdash t = s$ if $t^v = s^v$, and we extend this relation to all
formulas of $L_{\M}$ by the usual rules for the Boolean connectives. We set
$(\alg{X}, v) \Vdash \forall x\, \varphi$ if $(\alg{X}, v') \Vdash \varphi$
for every valuation $v'$ that agrees with $v$ outside $x$. A sentence
$\varphi$ is \emph{valid} on $\alg{X}$, in symbols $\alg{X} \Vdash \varphi$,
if $(\alg{X}, v) \Vdash \varphi$ for every valuation $v$. By construction,
the equivalence
\[
  \alg{X} \Vdash \varphi \quad\Longleftrightarrow\quad \alg{X}_+ \models
  \varphi
\]
holds for every sentence $\varphi$ of $L_{\M}$. We denote by
$\Fr(\varphi)$ the class of bounded poc sets on which a sentence $\varphi$ is
valid, whose members we call the \emph{frames} of $\varphi$.

The definition of $\Vdash$ can be written out in monadic second-order logic.
Let $L_{\Poc} := \{\leq, {}^\co, 0, 1\}$ be the first-order language of
bounded poc sets, whose element variables we denote by $p$ and $q$, possibly
indexed. Its monadic second-order extension $\mathrm{MSO}(L_{\Poc})$ has, in
addition, set variables $U$, which range over the subsets of the structure,
and atomic formulas $\tau \in U$ for every term $\tau$ of $L_{\Poc}$. The
\emph{standard translation} $\ST$, from the formulas of $L_{\M}$ to those of
$\mathrm{MSO}(L_{\Poc})$, is defined in the following steps.
\begin{enumerate}
\item With variables $x \neq y$ of $L_{\M}$ we associate set variables
  $U_x \neq U_y$.
\item With each term $t$ of $L_{\M}$ we associate a formula $\ST_p(t)$
  defined by induction on $t$. For a variable $x$ we set
  $\ST_p(x) := (p \in U_x)$, and for terms $t_1$, $t_2$ and $t_3$ we set
  \begin{multline*}
    \ST_p\bigl(m(t_1, t_2, t_3)\bigr) \\ := \bigl(\ST_p(t_1) \wedge
    \ST_p(t_2)\bigr) \vee \bigl(\ST_p(t_1) \wedge \ST_p(t_3)\bigr) \vee
    \bigl(\ST_p(t_2) \wedge \ST_p(t_3)\bigr).
  \end{multline*}
  Note that $p$ is a fixed element variable, the only free element variable of
  $\ST_p(t)$.
\item An atomic formula $t = s$ of $L_{\M}$ is translated as
  \[
    \ST(t = s) := \forall p\, \bigl(\ST_p(t) \leftrightarrow
    \ST_p(s)\bigr).
  \]
\item The translation commutes with the Boolean connectives.
\item Quantifiers over elements become quantifiers over valuations, that is,
  we set $\ST(\forall x\, \varphi) := \forall U_x\, (\mathrm{uf}(U_x) \to
  \ST(\varphi))$ and $\ST(\exists x\, \varphi) := \exists U_x\,
  (\mathrm{uf}(U_x) \wedge \ST(\varphi))$, where $\mathrm{uf}(U)$ is a
  formula of $\mathrm{MSO}(L_{\Poc})$ expressing that $U$ is an
  ultrafilter.\footnote{For instance, take $\mathrm{uf}(U) := \forall p
  \forall q\, ((p \in U \wedge p \leq q) \to q \in U) \wedge \forall p\,
  (p \in U \leftrightarrow \neg(p^\co \in U))$.}
\end{enumerate}
An induction on
formulas gives
\begin{equation}\label{eq:ST}
  (\alg{X}, v) \Vdash \varphi \quad\Longleftrightarrow\quad \alg{X} \models
  \ST(\varphi)\bigl(v(x_1), \ldots, v(x_n)\bigr)
\end{equation}
for every formula $\varphi(x_1, \ldots, x_n)$ of $L_{\M}$ and every valuation
$v$ on $\alg{X}$, each morphism $v(x_i)$ being identified with the ultrafilter
that it characterizes. In particular, a sentence $\varphi$ is valid on $\alg{X}$
if and only if $\alg{X} \models \ST(\varphi)$. Every sentence of $L_{\M}$ thus
corresponds, on bounded poc sets, to a sentence of monadic second-order logic
over $L_{\Poc}$. For some sentences more can be said, because the quantifiers
over valuations can be eliminated from their translation, which is then
equivalent to a first-order sentence of $L_{\Poc}$, that is, to an elementary
property of the poc set.

\begin{example}\label{ex:rho-k-correspondent}
The sentences $\rho_k$ of Example~\ref{ex:rho-k} with $k \geq 1$ have
first-order $L_{\Poc}$ correspondents. Two elements $p$ and $q$ of a bounded poc
set are \emph{transverse} if none of $p \leq q$, $p \leq q^\co$, $p^\co \leq q$
and $q \leq p$ holds. We denote by $\Trans(p, q)$ the quantifier-free formula of
$L_{\Poc}$ that expresses this condition. Applied to the poc set $X \setminus
\{0, 1\}$, in which transverse elements are called crossing,
\cite[Lemma~9.3.1]{Bowditch2022} states that the rank of $\alg{X}_+$ is the
largest number of pairwise transverse elements of $X \setminus \{0, 1\}$. By
\eqref{eq:ST}, for every bounded poc set $\alg{X}$,
\[
\alg{X} \models \ST(\rho_k) \quad\Longleftrightarrow\quad \alg{X} \models
\forall p_0 \cdots \forall p_k \bigvee_{0 \leq i < j \leq k}
\neg \Trans(p_i, p_j).
\]
The sentence $\ST(\rho_k)$ quantifies over $2^{k+1}$ morphisms into
$\utilde{\alg{2}}$, and its first-order correspondent over $k+1$ elements of
the poc set.
\end{example}

The correspondent of Example~\ref{ex:rho-k-correspondent} is an instance of a
general fact, which covers every universal sentence $\varphi$ with
$\classop{H}(\Mod(\varphi)) = \Mod(\varphi)$, as $\rho_k$ is
(Example~\ref{ex:rho-k}). The closure of $\Mod(\varphi)$ under
homomorphic images corresponds to the closure of $\Fr(\varphi)$ under
substructures.

\begin{proposition}\label{prop:universal-correspondent}
Let $\varphi = \forall x_1 \cdots \forall x_n\, \psi$ be a universal sentence of
$L_{\M}$, with $n \geq 1$ and $\psi$ quantifier-free.
\begin{enumerate}
\item\label{it:corr-small} If $\alg{X} \not\Vdash \varphi$ for a bounded poc
  set $\alg{X}$, then $\alg{Y} \not\Vdash \varphi$ for some substructure
  $\alg{Y}$ of $\alg{X}$ with at most $2^n$ elements.
\item\label{it:corr-dual} The equalities
$\classop{H}(\Mod(\varphi)) = \Mod(\varphi)$
and
$\classop{S}(\Fr(\varphi)) = \Fr(\varphi)$
  are equivalent.
\item\label{it:corr-closed} If these equalities hold, then $\varphi$ is valid
  on a bounded poc set $\alg{X}$ if and only if it is valid on every
  substructure of $\alg{X}$ with at most $2^n$ elements, and the translation
  $\ST(\varphi)$ is equivalent, over bounded poc sets, to a universal sentence
  of $L_{\Poc}$.
\end{enumerate}
\end{proposition}
\begin{proof}[Sketch of proof]
We view the elements of $\alg{X}_+$ as the morphisms $\alg{X} \to
\utilde{\alg{2}}$. An embedding $g \colon \alg{Y} \to \alg{X}$ then induces
the homomorphism $g_+ \colon u \mapsto u \circ g$ from $\alg{X}_+$ onto
$\alg{Y}_+$. For finite bounded poc sets, this is the injectivity of the alter
ego of the strong duality of Theorem~\ref{thm:dualities}(\ref{it:dual-poctau}),
and the general case follows from the fact that $\alg{X}$ is the directed
union of its finite substructures.

For~(\ref{it:corr-small}), let $\psi$ fail in $\alg{X}_+$ at $a_1, \ldots,
a_n$. Evaluation at the elements of $X$ yields at most $2^n$ homomorphisms from
the subalgebra $\alg{B}$ generated by $a_1, \ldots, a_n$ into $\alg{2}$, and a
substructure $\alg{Y}$ of $\alg{X}$ that contains one element for each of them
has at most $2^n$ elements and already separates the elements of $B$. Hence
$\alg{B}$ embeds into $\alg{Y}_+$, where $\psi$ fails. In~(\ref{it:corr-dual}),
the first equality implies the second because $\alg{Y}_+ \in
\classop{H}(\alg{X}_+)$ for every substructure $\alg{Y}$ of $\alg{X}$. The
converse follows by restricting a surjective homomorphism $f \colon \alg{A}
\to \alg{B}$ to finitely generated, hence finite, subalgebras, which are the
complex algebras of their poc sets, and by embedding $\alg{B}^+$ into
$\alg{A}^+$ through precomposition with $f$. Item~(\ref{it:corr-closed})
follows from~(\ref{it:corr-small}) and~(\ref{it:corr-dual}), since forbidding
the finitely many bounded poc sets with at most $2^n$ elements on which
$\varphi$ is not valid is expressed by a universal sentence of $L_{\Poc}$.
\end{proof}

A universal sentence of $L_{\Poc}$ that holds in a bounded poc set holds in
its substructures. By Proposition~\ref{prop:universal-correspondent}, the
translation of a universal sentence $\varphi$ of $L_{\M}$ is thus equivalent
to a universal sentence of $L_{\Poc}$ if and only if
$\classop{H}(\Mod(\varphi)) = \Mod(\varphi)$. The sentence
$\sigma$ of Remark~\ref{rem:transfer-not-necessary} shows that this equality
may fail, and $\ST(\sigma)$ is equivalent to no universal sentence of
$L_{\Poc}$.

\begin{corollary}\label{cor:positive-universal}
Let $\varphi$ be a universal sentence of $L_{\M}$ such that
$\classop{H}(\Mod(\varphi)) = \Mod(\varphi) \neq \emptyset$. Over median
algebras, the sentence $\varphi$ is equivalent to a positive universal
sentence.
\end{corollary}
\begin{proof}
Write $\varphi = \forall x_1 \cdots \forall x_n\, \psi$, with $\psi$
quantifier-free. By local finiteness, there are finitely many pairs
$(\alg{C}, \mathbf{c})$, up to isomorphism, in which $\mathbf{c}$ is an
$n$-tuple that generates $\alg{C}$ and at which $\psi$ fails, and, for the same
reason, finitely many terms in the variables $x_1, \ldots, x_n$ up to
equivalence over median algebras. The trivial algebra lies in $\Mod(\varphi)$,
as a homomorphic image of any of its members, so that each such $\alg{C}$ has
at least two elements. The positive universal sentence
\[
  \chi_{\alg{C}, \mathbf{c}} := \forall x_1 \cdots \forall x_n \bigvee \{s =
  t \mid s(\mathbf{c}) \neq t(\mathbf{c})\}
\]
fails at a tuple $\mathbf{a}$ exactly when $\mathbf{a} \mapsto \mathbf{c}$
extends to a homomorphism from the subalgebra generated by $\mathbf{a}$ onto
$\alg{C}$. Since $\Mod(\varphi)$ is closed under subalgebras and homomorphic
images, the sentence $\varphi$ is thus equivalent to the conjunction of the
sentences $\chi_{\alg{C}, \mathbf{c}}$.
\end{proof}

By Corollary~\ref{cor:positive-universal}, items~(\ref{it:transfer-closed})
and~(\ref{it:transfer-pos}) of Theorem~\ref{thm:transfer} apply to the
same classes of median algebras, item~(\ref{it:transfer-closed}) being the form
used in practice, as for the sentences $\rho_k$, which are not positive as
written. For universal sentences $\varphi$
with $\classop{H}(\Mod(\varphi)) = \Mod(\varphi)$,
Proposition~\ref{prop:universal-correspondent} and
Theorem~\ref{thm:transfer}(\ref{it:transfer-closed}) thus settle both the
existence of a universal correspondent and the transfer to the profinite
completion, so that a median algebra $\alg{A}$ satisfies $\varphi$ if and only
if $\hat{\alg{A}}$ does, that is, if and only if $\varphi$ is valid on every
substructure of $\alg{A}^+$ with at most $2^n$ elements.

What these results leave open is the shape of the correspondents. For the
sentences $\rho_k$, Proposition~\ref{prop:universal-correspondent} forbids
substructures with up to $2^{2^{k+1}}$ elements, whereas
Example~\ref{ex:rho-k-correspondent} only forbids $k+1$ pairwise transverse
elements. The proposition thus provides correspondents that are computable from
the sentence, but that are neither small nor readable from its syntax.

Since median algebras are locally finite, every universal sentence $\varphi$
in $n$ variables is equivalent over median algebras to the conjunction of the
sentences $\varphi_{\alg{B}}$ of Example~\ref{ex:rho-k}, where $\alg{B}$
ranges, up to isomorphism, over the finitely many median algebras generated by
at most $n$ elements in which $\varphi$ fails. For instance, the sentence
$\sigma$ of Remark~\ref{rem:transfer-not-necessary} is equivalent to
$\varphi_{\alg{S}}$ for the star $\alg{S}$ with four leaves. This suggests the
following questions.

\begin{question}\label{q:correspondence}
\begin{enumerate}
\item For which finite median algebras $\alg{B}$ does the equality
  $\classop{H}(\Mod(\varphi_{\alg{B}})) =
  \Mod(\varphi_{\alg{B}})$ hold? Is this the case exactly when the dual
  $\alg{B}^*$ is injective in $\Poc_\tau$, as it is for cubes
  (Lemma~\ref{lem:cube-lift})?
\item Is there a syntactically defined class of sentences of $L_{\M}$, beyond
  the universal ones, whose translations are equivalent to first-order
  sentences of $L_{\Poc}$ computed from their syntax? Is $\ST(\sigma)$
  equivalent to a first-order sentence of $L_{\Poc}$?
\item Does every universal sentence that holds in a median algebra $\alg{A}$
  hold in $\hat{\alg{A}}$, that is, does every finite subalgebra of
  $\hat{\alg{A}}$ embed into $\alg{A}$?
\item Is there a syntactically defined class of sentences of the form
  $\forall \mathbf{x}\, \exists \mathbf{y}\, \psi$, with $\psi$
  quantifier-free, that hold in $\hat{\alg{A}}$ whenever they hold in
  $\alg{A}$? Not all such sentences do.
\end{enumerate}
\end{question}

\bigskip

\paragraph{\bfseries Median algebra expansions} Equational, and even
quasi-equational, definability is degenerate for median algebras, since every
median algebra with at least two elements contains a copy of $\alg{2}$, so that
$\M$ is the only nontrivial quasivariety of median algebras. Universal
sentences, such as the sentences $\rho_k$, already define proper nontrivial
classes. Those universal sentences $\varphi$ with $\classop{H}(\Mod(\varphi)) =
\Mod(\varphi)$ are preserved by profinite completion
(Theorem~\ref{thm:transfer}), while the case of the other ones is the third
question of Question~\ref{q:correspondence}. Sentences with alternating
quantifiers need not be preserved either.

For expansions of median algebras by further operations, equational definability
becomes meaningful. When the new operations are homomorphisms, they extend to
the profinite completion by Proposition~\ref{prop:ump-reflection}, and the
definability questions can be asked as before. Retract algebras are an
interesting example. A \emph{retraction} of a median algebra is an idempotent
homomorphism with convex image, and a map $r$ is a retraction if and only if \[
r(m(x, y, z)) = m(r(x), r(y), z) \]
\cite[Proposition~5.1]{BandeltHedlikova1983}.
The algebras $(A, m, r)$ satisfying this equation form a variety, with a
nontrivial equational theory. For instance, the equation $r(x) = r(y)$ defines
the subvariety of constant retractions, which is term-equivalent to the pointed
median algebras of Definition~\ref{def:cat-ms},
while the equation $r(x) = x$ defines the subvariety of median algebras. A
retraction $r$ extends to a continuous homomorphism $\hat{r}$ of
$\hat{\alg{A}}$, and every variety of retract algebras is thus closed under
profinite completion. The dual map $r^+$ is an endomorphism of $\alg{A}^+$
that fixes the homomorphisms $u \colon \alg{A} \to \alg{2}$ that are not
constant on the image of $r$, and sends the other ones to the constant
homomorphisms. This leads us to the following question.

\begin{question}\label{q:retract}
Develop the correspondence theory of varieties of retract algebras, and a
Goldblatt--Thomason theorem for them. Their frames are the bounded poc sets
equipped with an idempotent morphism $e$ that sends every element either to
itself or to a bound.
\end{question}

Grau's ternary Boolean algebras \cite{Grau1947} are another instance, the
complement being an involutive automorphism of the median reduct. Its
continuous extension to $\hat{\alg{A}}$ satisfies Grau's equations, which hold
on the dense subalgebra $\alg{A}$, so that the profinite completion of the
median reduct of a ternary Boolean algebra, equipped with the extended
complement, is again a ternary Boolean algebra.

Operations that are not homomorphisms raise a problem that does not occur for
retractions, namely the choice of their extension to $\hat{\alg{A}}$. Every map
between median algebras has a canonical multi-valued extension
\cite[Definition~4.2]{HansoulTeheux}, which the conditional completeness of
$(\hat{\alg{A}}, \leq_a)$ (Proposition~\ref{prop:cond-complete}) turns into a
single-valued one by taking meets, in the spirit of the $\sigma$-extensions of
lattice expansions \cite{Gehrke2004}. The resulting extension depends however on
the basepoint $a$, unlike in the lattice setting.

\begin{question}\label{q:expansions}
Which equations of the median algebra expansions are preserved by profinite
completion, for instance for expansions by operations that preserve the joins in
each interval of $(\alg{A}, \leq_a)$? When does the extension of an operation
not depend on the basepoint?
\end{question}

Pointed median algebras build the dependence on the basepoint into the language.
By Proposition~\ref{prop:ms-pointed}, they amount to median semilattices $(A,
\wedge_a)$, in which the meet is a total operation while joins exist only for
pairs that are bounded above. We therefore call a \emph{modal operator} a map
$\Box$ satisfying $\Box(x \wedge_a y) = \Box x \wedge_a \Box y$. As there is in
general no top element to preserve, a modal operator is a regular rather than a
normal box. Retractions provide examples. If $a$ lies in the image of a
retraction $r$, then $r$ preserves the meets $m(a, x, y)$, being a homomorphism
that fixes $a$. Moreover $r(x) \leq_a x$ and $r(r(x)) = r(x)$, so that $r$ is an
interior operator, as the box of the modal logic $\mathsf{S4}$. Every interval
of $(A, \leq_a)$ being a distributive lattice, Theorem~\ref{thm:can-ext}
suggests asking canonicity interval by interval, as for distributive lattices
with operators \cite{Gehrke2004}, and the conditional canonical extension
$\hat{\alg{A}}_a$ of Section~\ref{sec:ccan} is a natural candidate for the
completion.

The dual notion is a map $\Diamond$ that preserves the bottom $a$ and the
existing joins of $(A, \leq_a)$. In the absence of a complement, it cannot be
defined from a modal operator, and as joins are partial, its defining
property is not an equation. By Lemma~\ref{lem:interval-dual}, however, the
joins for $\leq_a$ in an interval $[a; d]$ are the meets for $\leq_d$. As
$\Diamond$ fixes $a$ and preserves $\leq_a$, it maps $[a; d]$ into $[a;
\Diamond d]$, and it turns the meets of $([a; d], \leq_d)$ into those of $([a;
\Diamond d], \leq_{\Diamond d})$. The duality between $\Diamond$ and $\Box$
thus passes through a change of basepoint, from $d$ to $\Diamond d$, interval
by interval.
\begin{question}\label{q:modal-ms}
\begin{enumerate}
\item Is there a Sahlqvist-type class of equations in $\wedge_a$ and a modal
  operator $\Box$, preserved when $\Box$ is extended to $\hat{\alg{A}}_a$ or
  to $\hat{\alg{A}}$, that contains the equations satisfied by retractions?
\item Let $\Diamond$ be a map that preserves the bottom $a$ and the existing
  joins of $(A, \leq_a)$. How does it extend to the completions, and how does
  the canonicity of its equations relate to that of modal operators for the
  other basepoints?
\end{enumerate}
\end{question}

\bigskip

\paragraph{\bfseries Cuboid topology}
In the polyhedral semantics of modal logic, propositions are interpreted as
subpolyhedra of a polyhedron and the box as the topological interior
\cite{BezhanishviliMarraMcNeillPedrini2018}.
A $\mathrm{CAT}(0)$ cube complex is a polyhedron whose cubes correspond, in
the median algebra $\alg{A}$ of its vertices, to the intervals isomorphic to
some $\alg{2}^n$. For a basepoint $a$, the map $x \mapsto [a;x]$ sends the
meet $\wedge_a$ to the intersection of convex subcomplexes, which are unions
of cubes. This suggests using modal expansions of pointed median algebras as
languages for the topology of cube complexes. If the connection between such
languages and this cuboid topology can be made precise, for instance by
relating the rank conditions $\rho_k$ to dimension, it would open a research
program of its own.

\section*{Declaration of generative AI use}

During the preparation of this work, the author used generative AI models to
assist with drafting and editing text, checking arguments, searching the
literature (in particular that of geometric group theory), finding some of
the counterexamples, and investigating some of the questions raised in
Section~\ref{sec:conclusion}. All results and statements were checked by the
author, who takes full responsibility for the content of the paper.

\bibliographystyle{amsalpha}
\bibliography{Notes_definability_pocset}

\end{document}